\documentclass[11pt]{article}
\usepackage{amsmath,amssymb,latexsym}

\newtheorem{thm}{Theorem}
\newtheorem{lem}[thm]{Lemma}
\newtheorem{cor}[thm]{Corollary}
\newtheorem{rmk}[thm]{Remark}

\newtheorem{dfn}[thm]{Definition}
\newtheorem{exm}{Example}
\newcommand{\bea}{\begin{eqnarray}}
\newcommand{\eea}{\end{eqnarray}}
\newcommand{\bean}{\begin{eqnarray*}}
\newcommand{\eean}{\end{eqnarray*}}
\newcommand{\beq}{\begin{equation}}
\newcommand{\eeq}{\end{equation}}
\newcommand{\bac}{\begin{array}{c}}
\newcommand{\ball}{\begin{array}{ll}}
\newcommand{\ea}{\end{array}}

\def\({\left(}
\def\){\right)}

\title{On the Whitney Extension-Interpolation-Alignment problem for almost isometries with small distortion in $\Bbb R^D$}
\author{Steven B. Damelin \thanks{Zentralblatt MATH, FIZ Karlsruhe – Leibniz Institute\, email: steve.damelin@gmail.com} \and Charles Fefferman, \thanks{Department of Mathematics; Fine Hall, Washington Road, Princeton NJ 08544-1000 USA,\, email: cf@math.princeton.edu}}
\date{}

\begin{document}
\maketitle

\thispagestyle{empty}
\parskip=10pt

\begin{abstract}
\noindent

Let $D\geq 2$, $S\subset \mathbb R^D$ be finite and let $\phi:S\to \mathbb R^D$ with
$\phi$ a small distortion on $S$. We solve the Whitney extension-interpolation-alignment problem of how to understand when $\phi$ can be extended to a function 
$\Phi:\mathbb R^D\to \mathbb R^D$ which is a smooth small distortion on $\mathbb R^D$. Our main results are in addition to 
Whitney extensions, results on interpolation and alignment of data in $\mathbb R^D$ and complement those of \cite{FD3}. The work in this paper appears in the research memoir \cite{SDam}.
\end{abstract}
\medskip

\noindent AMS-MSC Classification: 58C25, 42B35, 94A08, 94C30, 41A05, 68Q25, 30E05, 26E10, 68Q17, 
\vskip 0.2in
\tableofcontents

\section{Whitney's extension problem.}
\setcounter{equation}{0}
We will work here and throughout in Euclidean space $\mathbb R^D$, $D\geq 2$ and $\mathbb R^n$, $n\geq 1$. $|.|$ will denote the Euclidean norm on $\mathbb R^D$ or $\mathbb R^n$ if not stated otherwise.
Henceforth, $D,n$ are chosen and fixed. Throughout $K\geq 1$ is a special constant, chosen and fixed, depending only on $D$ and will be defined in (\ref{e:specialk}). 
\medskip

{\it Question 1}\, Let $m\geq 1$ and let $\phi:S\to \mathbb R$ be a function defined on an arbitrary set $S\subset \mathbb R^n$. How can one decide whether $\phi$ extends to a function $\Phi:\mathbb R^n\to \mathbb R$ which agrees with $\phi$ on $S$ and is in $C^m(\mathbb R^D)$, the space of functions from $\mathbb R^n$ to $\mathbb R$ whose derivatives of order $m$ are continuous and bounded.

For $n=1$ and $S$ compact, this is the well-known Whitney extension problem. See \cite{W,W1,W2}. Continued progress on this problem was made by G. Glaeser \cite{Gl}, Y. Brudnyi and P. Shvartsman \cite{B,BS1,BS2,BS3,BS4} and E. Bierstone, P. Milman and W. Pawlucki \cite{BMP}. See 
also N. Zobin \cite{Z,Z1}, B. Klartag and N. Zobin \cite{KZ} and 
E. LeGryer \cite{Le} for related work. C. Fefferman and his collaborators, A. Israel, B. Klartag, G. Lui, S. Mitter and H. Narayanan \cite{F1,F3,F4,F5,F6,KF,KF1,KF2,KF3, FIL1,FIL2,FIL3,FIL4,FIL5,FMH,I} have given a complete solution to this problem and have extended and generalized it in many ways. 

\subsection{Whitney's extension problem for small distortions in $\mathbb R^D$.}
In this paper we solve the following Whitney extension problem, continuing from \cite{FD3} given in:

{\it Question 2}\, Let $S\subset \mathbb R^D$ be a finite set with some geometry to be defined more precisely later. Suppose that we are given a map $\phi:S\to \mathbb R^D$ with
$\phi$ a small distortion on $S$. How can one decide whether $\phi$ extends to a smooth small distortion $\Phi:\mathbb R^D\to \mathbb R^D$ which agrees with $\phi$ on $S$. Our main results are in addition to Whitney extensions, results on interpolation and alignment of data in $\mathbb R^D$. 

\section{Main Results: Whitney extensions-Interpolation-Alignment in $\mathbb R^D$.}
\setcounter{equation}{0}
\subsection{Three needed definitions.}
Our answer to Question 2 comes in the form of our five  main results below which rely on three crucial objects. Throughout, given a positive number $\varepsilon$, 
we assume that $\varepsilon$ is less than a small enough constant determined by $D$ and $K$, see ((\ref{e:specialk})). We call this the small $\varepsilon$ assumption. Given a positive number $\varepsilon'$, we assume that $\varepsilon'$ is less than a small enough constant determined by $n$. We call this the small $\varepsilon'$ assumption.

Our first object is needed to describe a class of smooth 1-1 and onto small distortions from $\mathbb R^D \to\mathbb R^D$ ($\mathbb R^n \to\mathbb R^n$). 
which we call $\varepsilon$-distorted diffeomorphisms ($\varepsilon'$-distorted diffeomorphisms), the later first introduced in \cite{FD3}. The second object is needed to describe the geometry of $S$. We call it an $\eta$ block. The third is a Euclidean motion $A: \mathbb R^n \to \mathbb R^n$. 
\begin{dfn}
A diffeomorphism (hence 1-1 and onto) $\Phi:\mathbb R^n\to \mathbb R^n$ is ``$\varepsilon$' distorted" provided
\beq
(1+\varepsilon')^{-1}I\leq [\nabla \Phi(x)]^{T}[\nabla \Phi (x)]\leq (1+\varepsilon')I
\label{e:edistortion}
\eeq
as matrices, for all $x\in \mathbb R^n$. Here, $I$ denotes the identity matrix in $\mathbb R^n$. 
\end{dfn}

If $\Phi$ is $\varepsilon'$ distorted then an application of Bochner's theorem gives that we have for all $x,y\in \mathbb R^n$,
\beq
(1+\varepsilon')^{-1}|x-y|\leq |\Phi(x)-\Phi(y)|\leq |x-y|(1+\varepsilon').
\label{e:edistortion1}
\eeq

A $\varepsilon'$-distorted diffeomorphism $\Phi$ is proper if ${\rm det}(\Phi')>0$ on $\mathbb R^n$ and improper if ${\rm det}(\Phi')<0$
on $\mathbb R^n$. Since ${\rm det}(\Phi')\neq 0$ everywhere on $\mathbb R^n$, every $\varepsilon'$- distorted
diffeomorphism is proper or improper. $\varepsilon$-distorted diffeomorphisms are defined similarly. 

Henceforth, by a Euclidean motion on $\mathbb R^n$, we shall mean a map $A:x\to Ex+x_0$ from $\mathbb R^n$ to $\mathbb R^n$ with $E\in O(n)$ or $E\in SO(n)$ and $x_0\in \Bbb R^n$ fixed. $O(n)$ and $SO(n)$ will denote the orthogonal and special orthogonal groups respectively. A Euclidean motion $A$ will be called "proper" if $E\in SO(n)$ otherwise it is called improper. 
\footnote{$O(n)$ is by definition the group of isometries of $\mathbb R^n$ which preserve a fix point and $SO(n)$ the subgroup of $O(n)$ of orthogonal matrices of determinant 1. The Euclidean motions of $\mathbb R^n$ are the elements of the symmetry group of $\mathbb R^n$, ie all isometries of $\mathbb R^n$.}

More generally, an invertible affine map $T:\mathbb R^D\to \mathbb R^D$ is proper if ${\rm det}(T')>0$ and improper if ${\rm det}(T')<0$. Since $T$ is invertible, $T$ is either proper or improper.

\begin{dfn}
{\rm For $z_0,z_1,...,z_l \in \mathbb R^D$ with $l\leq D$, $V_l(z_0,...,z_l)$ will denote the $l$-dimensional volume of the $l$-simplex with vertices at $z_0,...,z_l$. If $S\subset \mathbb R^D$ is a finite set, then $V_l(S)$ denotes the max of $V_l(z_0,...,z_l)$ over all 
$z_0,z_1,...,z_l\in S$. \footnote{If $V_D(S)$ is small, then we expect that $S$ will be close to a hyperplane in $\mathbb R^D$.} Let 
$\phi:S\to \mathbb R^D$ and let $0<\eta<1$. A positive (resp. negative) $\eta$-block for $\phi$ is a $D+1$ tuple $(x_0,...,x_D)\in \mathbb R^D$ such that the following two conditions hold: (1) $V_D(x_0,...,x_D)\geq (\leq) \eta^D{\rm diam}(x_0,...,x_D)$. (2) Let $T$ be the unique affine map which agrees with $\phi$ on $S$. $T$ exists and is unique by virtue of \cite{ATV}. Then we assume that $T$ is proper or improper. (Note that if the map $T$ is not invertible then $(x_0,...,x_D)$ is not an 
$\eta$ block.) }
\label{d:block}
\end{dfn}

\subsection{Five main results: Whitney extensions-interpolation and alignment in $\mathbb R^D$.}
Our five main results of this paper which provide an answer to Question 2 are as follows:

\begin{thm}
Let $S\subset \mathbb R^D$ be finite. There exists positive constants $c_K$, $C'_K$, $C''_K$ depending only on $D$ and $K$ such that the 
following holds: Set $\eta=\exp(-C'_K/\varepsilon)$ and $\delta=\exp(-C''_K/\varepsilon)$ with $0<\varepsilon<c_K$.  Let $\phi:S\to \mathbb R^D$ satisfy
\beq
(1+\delta)^{-1}|x-y|\leq |\phi(x)-\phi(y)|\leq (1+\delta)|x-y|, \, x, y\in S.
\label{e:deltadistortion}
\eeq
Then if $\phi$ has no negative $\eta$ block, there exists a proper $\varepsilon$-distorted diffeomorphism
$\Phi:\mathbb R^D\to \mathbb R^D$ such that $\phi=\Phi$ on $S$ and $\Phi$ agrees with a proper Euclidean motion on
\[
\left\{x\in \mathbb R^D:\, {\rm dist}(x,S)\geq 10^4{\rm diam}(S)\right\}.
\]
\label{t:Theorem2a}
\end{thm}

\begin{thm}
Let $S\subset \mathbb R^D$ be finite. There exists positive constants $c_K$, $C'_K$, $C''_K$ depending only on $D$ and $K$ such that the 
following holds: Set $\eta=\exp(-C'_K/\varepsilon)$ and $\delta=\exp(-C''_K/\varepsilon)$ with $0<\varepsilon<c_K$. Let $\phi:S\to \mathbb R^D$ satisfy
$(\ref{e:deltadistortion})$. Then if $\phi$ has a negative $\eta$ block, $\phi$ cannot be extended to a proper $\delta$ distorted diffeomorphism of $\mathbb R^D$.
\label{t:Theorem2b}
\end{thm}

\begin{thm}
Let  $S\subset \mathbb R^D$ be finite. There exists positive constants $c_K$, $C'_K$ depending only on $D$ and $K$ such that the 
following holds: Set $\delta=\exp(-C'_K/\varepsilon)$ with $0<\varepsilon<c_K$.  Let $\phi:S\to \mathbb R^D$ satisfy
$(\ref{e:deltadistortion})$. Suppose that for any $S_o\subset S$ with at most $2D+2$ points, there exists a $\delta$ distorted diffeomorphism $\Phi^{S_0}:\mathbb R^D\to \mathbb R^D$ such that $\Phi^{S_0}=\phi$ on $S_0$. Then, there exists an $\varepsilon$-distorted diffeomorphism $\Phi:\mathbb R^D\to \mathbb R^D$ such that $\Phi=\phi$ on $S$.
\label{t:Theorem3}
\end{thm}

\begin{thm}
Let $S\subset \mathbb R^D$ with ${\rm card}(S)\leq D+1$. There exist positive constants $c,C$ depending only on $D$ such that the following holds: Set 
$\delta=\exp(-C/\varepsilon)$ with $0<\varepsilon<c$ and let $\phi:S\to \mathbb R^D$ satisfy
$(\ref{e:deltadistortion})$.
Then there exists a $\varepsilon$-distorted diffeomorphism $\Phi:\mathbb R^D\to \mathbb R^D$ such that $\Phi=\phi$ on $S$.
\label{t:Theorem4}
\end{thm}

\begin{thm}
Let $\phi:E\to \mathbb R^D$ where $E\in \mathbb R^D$ is finite and let $0<\eta<1$. Suppose that $\phi$ satisfies $(\ref{e:deltadistortion})$ and has a positive $\eta$ block and a negative
$\eta$ block. Let $0<\delta<c\eta^D$ for small enough $c>0$ depending only on $D$. Then $\phi$ does not extend to a
$\delta$ distorted diffeomorphism $\Phi:\mathbb R^D\to \mathbb R^D$.
\label{t:cextensionblock2}
\end{thm}

\section{Discussion of Main Results.}
\setcounter{equation}{0}
This section is devoted to understanding our main results, Theorem~\ref{t:Theorem2a}, Theorem~\ref{t:Theorem2b}, Theorem~\ref{t:Theorem3},
Theorem~\ref{t:Theorem4} and Theorem~\ref{t:cextensionblock2}.

\subsection{Comparisons. }
In \cite{FD3}, we proved the following Whitney Extension Theorem:

\begin{thm}
Let $\varepsilon'>0$ and let $y_1,...y_k$ and $z_1,...z_k$ be two $1\leq k\leq n$ sets of distinct points in $\mathbb R^n$. Then there exists $\delta'>0$ depending only on $\varepsilon'$ such that the following holds: Suppose that
\beq
(1+\delta')^{-1}\leq \frac{|z_i-z_j|}{|y_i-y_j|}\leq (1+\delta'),\, 1\leq i,j\leq k,\, i\neq j.
\label{e:emotionsa}
\eeq
Then there exists a $\varepsilon'$-distorted diffeomorphism $\Phi:\mathbb R^n\to \mathbb R^n$ satisfying
\beq
\Phi(y_i)=z_i,\, 1\leq i\leq k.
\label{e:emotionsaaa}
\eeq
\label{t:lemmafive}
\end{thm}

\begin{itemize}
\item[(a)] Note that Theorem~\ref{t:lemmafive} works for any two distinct $1\leq k\leq n$ configurations $y_1,...,y_k$, $z_1,...,z_k$ satisfying the condition (\ref{e:emotionsa}). If we assume that there exists a map $\phi:S\to \mathbb R^n$ 
where $S:=\left\{y_1,y_2,...,y_k\right\}$ and $\phi(S):=\left\{z_1,z_2,...,z_k\right\}$, then Theorem~\ref{t:lemmafive} is a {\bf Whitney extension theorem in $\mathbb R^n$}. In the same spirit, Theorem~\ref{t:Theorem2a}, Theorem~\ref{t:Theorem2b}, Theorem~\ref{t:Theorem3}, 
Theorem~\ref{t:Theorem4} and Theorem~\ref{t:cextensionblock2} are Whitney extension results in $\mathbb R^D$.
\item[(b)] Notice that Theorem~\ref{t:lemmafive} is an {\bf Interpolation of data result in $\mathbb R^n$} with the map $\Phi$ the interpolation operator. 
In the same spirit, Theorem~\ref{t:Theorem2a}, Theorem~\ref{t:Theorem2b}, Theorem~\ref{t:Theorem3}, 
Theorem~\ref{t:Theorem4} and Theorem~\ref{t:cextensionblock2} are interpolation results in $\mathbb R^D$.
If $z_k=\phi(y_k),k\geq 1$ with a $C^{m}(U), m\geq 1$, $\phi:U\to \mathbb R^D$ (some $U\subseteq \mathbb R^D$ open), whose restriction to $S$:=$(\left\{y_1,...y_k\right\}$ with some geometry) satisfies that $\frac{|\phi(y_i)-\phi(y_j)|}{|y_i-y_j|}$ is close to 1 for all $i\neq j$,
then the $C^{m}(U)$ norm of $\Phi$ when $\Phi$ exists, should not be too large in order to approximate $\phi$ by $\Phi$ on subsets of $U$ if possible. The size of the $C^{m}(U)$ norm of $\Phi$ is addressed by us in a future paper. 
We refer the reader to the papers by Lubinsky and his collaborators and the many references cited 
therein for a good insight into the subject of interpolation by polynomials \cite{L1,L2,L3,L4,L5,L6,L7,DL,DL1}. The size and optimal size of the interpolation/extension norm and the operator which achieves these (very often not polynomial) in Whitney extension problems is a difficult problem in general and has 
been addressed substantially by Fefferman and his collaborators in \cite{F1,F3,F4,F5,F6,KF,KF1,KF2,KF3, FIL1,FIL2,FIL3,FIL4,FIL5,FMH,I} and several other papers referenced therein.
Theorem~\ref{t:lemmafive}, Theorem~\ref{t:Theorem2a}, Theorem~\ref{t:Theorem2b}, Theorem~\ref{t:Theorem3}, 
Theorem~\ref{t:Theorem4} and Theorem~\ref{t:cextensionblock2} are generalizations, in particular, of a well known result of John on extensions of isometries. See for example \cite{WW1}.
\end{itemize}

\subsection{The intriguing restriction in Theorem~\ref{t:lemmafive} on the number of points being bounded by $n$, the dimension of the Euclidean space $\mathbb R^n$.}
Theorem~\ref{t:lemmafive} has an intriguing feature, namely the restriction of the number of points $k$ to be bounded by the dimension of $\mathbb R^n$, $n$. Indeed, we showed in \cite{FD3} that Theorem~\ref{t:lemmafive} is false if $k>n$ by constructing a counter example which amounted to showing that the 
mapping swapping two given real numbers $\zeta$ and $-\zeta$ and fixing the number $1$ cannot be extended to a continuous bijection of the line.
Theorem~\ref{t:Theorem2a}, Theorem~\ref{t:Theorem2b}, Theorem~\ref{t:Theorem3}, 
Theorem~\ref{t:Theorem4} and Theorem~\ref{t:cextensionblock2} tell us that we may remove the restriction on $k$ if 
roughly we require that on any $D+1$ of the $k$ points which form a relatively voluminous simplex, the extension $\Phi$ is orientation preserving.

\subsection{The constants $\varepsilon$ and $\delta$ and how they relate to each other.}
Notice that Theorem~\ref{t:lemmafive} has two positive constants $\varepsilon'$ and $\delta'$. $\varepsilon'$ determines the distortion of the map $\Phi$ whereas $\delta'$ determines the distortion of the pairwise distances between point sets $y_1,...,y_k$ and $z_1,...,z_k$ and depends only on $\varepsilon'$.
When looking at Theorem~\ref{t:lemmafive}, an immediate question which comes to mind is what is the quantative dependence of $\delta'$ on $\varepsilon'$. We have:

\begin{thm}
Let $\varepsilon'>0$ and let $y_1,...y_k$ and $z_1,...z_k$ be two $1\leq k\leq n$ sets of distinct points in $\mathbb R^n$ with
\beq
\sum_{i\neq j}|y_i-y_j|^2+\sum_{i\neq j}|z_i-z_j|^2=1,\, y_1=z_1=0.
\label{e:lemmafourpoints}
\eeq
Then there exist constants $J,J'>0$ depending only on $n$ such the following holds: Set $\delta'=J'\varepsilon'^{J}$ and suppose 
\beq
||z_i-z_j|-|y_i-y_j||<\delta', \, i\neq j.
\label{e:lemmafourpointsa}
\eeq
Then there exists a $\varepsilon'$-distorted diffeomorphism $\Phi:\mathbb R^n\to \mathbb R^n$ satisfying
\beq
\Phi(y_i)=z_i,\, 1\leq i\leq k.
\label{e:emotionsaaa}
\eeq
\label{t:lemmafiveprime}
\end{thm}

Theorem~\ref{t:Theorem2a} is our main result and Theorem~\ref{t:Theorem2b}, Theorem~\ref{t:Theorem3}, 
Theorem~\ref{t:Theorem4} and Theorem~\ref{t:cextensionblock2} follow from it and its proof. These results show that one can take typically $\delta=\exp(-C_K/\varepsilon)$ with a positive constant $C_K$ depending only on $D$ and $K$. Theorem~\ref{t:Theorem3} holds with $\delta=C_K\varepsilon$ and this follows from the main result in \cite{FD1}. Clearly  $\delta=\exp(-C_K/\varepsilon)$ is a more refined estimate than $\delta'=J'\varepsilon'^{J}$ (ignoring the constants $C_K, J, J'$).  $\delta=C_K\varepsilon$  is clearly optimal. 

We have explained how  Theorem~\ref{t:lemmafive}, Theorem~\ref{t:lemmafiveprime}, Theorem~\ref{t:Theorem2a}, Theorem~\ref{t:Theorem2b}, Theorem~\ref{t:Theorem3}, 
Theorem~\ref{t:Theorem4} and Theorem~\ref{t:cextensionblock2} are Whitney extension and interpolation results in Euclidean space. We now show they are results on alignment of data in Euclidean space as well. To understand this, we begin with:

\subsection{Alignment of data in Euclidean space and the Procrustes problem.}
An important problem in computer vision is comparing $k\geq 1$ point configurations in $\Bbb R^n$. \footnote{In computer vision, the phrase "$k\geq 1$ point configuration" means $k\geq 1$ distinct points.} One way to think of this is as follows: Given two $k$ configurations in $\Bbb R^n$, do there exist combinations of rotations, translations, reflections and compositions of these which map the one configuration onto the other. This is the shape registration problem. \footnote{In computer vision, two point clouds are said to have the same ``shape" if there exist combinations of rotations, translations, reflections and compositions of these which map the one set of points onto the other. This typically is called the registration problem.} A typical application of this problem arises in image processing and surface alignment where it is often necessary to align an unidentified image to an image from a given data base of images for example parts of the human body in face or fingerprint recognition. Thus the idea is to recognize points (often called landmarks) by verifying whether they align to points in the given data base. 
An image in $\Bbb R^{n}$ does not change under Euclidean motions. 
Motivated by this, in this paper, we will think of shape preservation in terms of whether there exists a Euclidean motion which maps one $k$ point configuration onto a second. \footnote{The shape identification problem can be stated more generally for other transformation groups. 
We restrict ourselves in this paper to Euclidean motions.}
In the case of labelled data (where the data points in each set are indexed by the same index set), an old approach called the Procrustes approach \cite{G,G1} analytically determines a Euclidean motion which maps the first configuration close to the other (in a $L^2$) sense.
There are a variety of ways to label points. See for example \cite{LW,WW}. 
In terms of good algorithms to do alignment of this kind, 
one method of Iterative Closet Point (ICP) for example is very popular and analytically computes an optimal rotation for alignment of $k$ point configurations, see for example \cite{DL} and the references cited therein for a perspective. 
Researchers in geometric processing think of the problem of comparing point clouds or finding a distance between them as to asking how to deform one point cloud into the other (each point cloud represented by a collection of meshes) in the sense of saying they have the same shape. 
We refer the reader to the references \cite{P,LW,WW,BCSZ,BSAB,WSW,WS,ZS,SW,DL,LAD,LD1,LD2, L, L1, BLCPFPJD, L2, SWK, SSK, VLBRC, LRS} and the many references cited therein for a broad perspective on applications of this problem.
\medskip

One way to dig deeper into the Procrustes problem is to compare pairwise distances between labelled points. In this regard, the following result is well known. See for example \cite{ATV,WW1}.
\medskip

\begin{thm}
Let $y_1,...,y_k$ and $z_1,...,z_k$ be two $k\geq 1$ point configurations in $\mathbb R^n$. Suppose that
\[
|z_i-z_j|=|y_i-y_j|,\, 1\leq i,j\leq k,\, i\neq j.
\]
Then there exists a Euclidean motion $A:\Bbb R^n\to \Bbb R^n$ such that $A(y_i)=z_i,\, i=1,...,k.$ If $k\leq n$, then $A$ can be taken as proper.
\label{t:Theorem 1}
\end{thm}

Alignment of point configurations from their pairwise distances are encountered for example in X-ray crystallography and in the mapping of restriction sites of DNA. See \cite{P,SSL} and the references cited therein. (In the case of one dimension, this problem is known as the turnpike problem or in molecular biology, it is known as the partial digest problem). See the work of \cite{RS,LW} for example which deals with algorithms and their running time for such alignments. We mention that a difficulty in trying to match point configurations is the absence of labels in the sense that often one does not know which point to map to which. We will not deal with the unlabeled problem in this paper but see \cite{NDS}. 

\subsection{ Labelled approximate alignment in $\mathbb R^n$ in the case where pairwise distances are distorted.}
To study the labelled alignment data problem in the case where pairwise distances are distorted, we proved in \cite{FD3} an analogy of Theorem~\ref{t:Theorem 1} namely Theorem~\ref{t:lemmafour}:
\medskip

\begin{thm}\begin{itemize}\item[(a)] Given $\varepsilon'>0$, there exists $\delta'>0$ depending only on $\varepsilon'$, such that the following holds. Let $y_1,...,y_k$ and
$z_1,...,z_k$ be two $k\geq 1$ point configurations in $\mathbb R^n$ satisfying $(\ref{e:emotionsa})$.
Then, there exists a Euclidean motion $A:\mathbb R^n\to \mathbb R^n$ such that
\beq
|z_i-A(y_i)|\leq \varepsilon' {\rm diam}\left\{y_1,...,y_k\right\}
\label{e:emotionsb}
\eeq
for each $1\leq i\leq k$. If $k\leq n$, then we can take $A$ to be a proper. 
\item[(b)] Suppose now that $y_1,...,y_k$ and
$z_1,...,z_k$ are two $k\geq 1$ point configurations in $\mathbb R^n$ so that $(\ref{e:lemmafourpoints})$ holds.
Then there exist constants $J,J'>0$ depending only on $n$ such the following holds: Set $\delta'=J'\varepsilon'^{J}$ and suppose $(\ref{e:lemmafourpointsa})$.
Then, there exists a Euclidean motion $A:\mathbb R^n\to \mathbb R^n$ such that
\beq
|z_i-A(y_i)|\leq \varepsilon'. 
\label{e:emotionsu}
\eeq
\label{t:lemmafour}
\end{itemize}
\end{thm}

Notice that if we require each set of points $y_1,...,y_k$ and $z_1,...,z_k$ to be in a bounded set of controlled radius, then we are able to be specific about the relationship between $\varepsilon'$ and $\delta'$ in Theorem~\ref{t:lemmafour},
namely $\delta'=J'\varepsilon^{J}$. The proof of Theorem~\ref{t:lemmafiveprime} uses a careful reworking of Theorem~\ref{t:lemmafive} in \cite{FD3} using Theorem~\ref{t:lemmafour}, part (b). We do not provide the details here but we do provide the proof of Theorem~\ref{t:lemmafour} from \cite{FD3} in order to illustrate in particular, regarding (b), the use of an inequality called Lojasiewicz's inequality which we will need later.

\subsection{ Lojasiewicz's inequality and Proof of Theorem~\ref{t:lemmafour}.}
{\bf Proof of Theorem~\ref{t:lemmafour}} \ We first prove (a): Suppose not. Then for each $l\geq 1$, we can find points $y_1^{(l)},...,y_k^{(l)}$ and $z_1^{(l)},...,z_k^{(l)}$ in
$\mathbb R^D$ satisfying (\ref{e:emotionsa}) with $\delta=1/l$ but not satisfying (\ref{e:emotionsb}). Without loss of generality, we may suppose that ${\rm diam}\left\{y_1^{(l)},...,y_k^{(l)}\right\}=1$ for each $l$ and that $y_1^{(l)}=0$ and
$z_l^{(1)}=0$ for each $l$. Thus $|y_i^{(l)}|\leq 1$ for all $i$ and $l$ and
\[
(1+1/l)^{-1}\leq \frac{|z_i^{(l)}-z_j^{(l)}|}{|y_i^{(l)}-y_j^{(l)}|}\leq (1+1/l)
\]
for $i\neq j$ and any $l$.
However, for each $l$, there does not exist an Euclidean motion
$A$ such that
\beq
|z_i^{(l)}-A(y_i^{(l)})|\leq \varepsilon'
\label{e:emotionsc}
\eeq
for each $i$. Passing to a subsequence, $l_1,l_2,l_3,...,$ we may assume
\[
y_i^{(l_{\mu})}\to y_i^{\infty},\, \mu\to \infty
\]
and
\[
z_i^{(l_{\mu})}\to z_i^{\infty},\, \mu\to \infty.
\]
Here, the points $y_i^{\infty}$ and $z_i^{\infty}$ satisfy
\[
|z_i^{\infty}-z_j^{\infty}|=|y_i^{\infty}-y_j^{\infty}|
\]
for $i\neq j$. Hence, by Theorem~\ref{t:Theorem 1}, there is an Euclidean motion $A^*:\mathbb R^n\to \mathbb R^n$ such that $A*(y_i^{\infty})=z_i^{\infty}$. Consequently,
for $\mu$ large enough, (\ref{e:emotionsc}) holds with $l_{\mu}$. This contradicts the fact that for each $l$, there does not exist a $A$ satisfying (\ref{e:emotionsc}) with $l$.
Thus, we have proved all the assertions of (a) except that we can take $A$ to be proper if $k\leq n$. To see this, suppose that $k\leq n$ and let $A$ be an improper Euclidean motion such that
\[
|z_i-A(y_i)|\leq \varepsilon'{\rm diam}\left\{y_1,...,y_k\right\}
\]
for each $i$. Then, there exists an improper Euclidean motion $A^*$ that fixes $y_1,...,y_k$ and in place of $A$, we may use the map $A^* o A$ in (a) so (a) is proved. We now prove (b). Here we use an inequality
called Lojasiewicz's inequality, \cite{SJS} which allows us to control the upper bound estimate in (\ref{e:emotionsb}) and replace it by the upper bound in (\ref{e:emotionsu}) provided the points $y_1,...,y_k$ and $z_1,...,z_k$ each lie in bounded sets of controlled radius.
The Lojasiewicz's inequality says the following: Let $f:U\to \Bbb R$ be a real analytic function on an open set $U$ in $\Bbb R^n$ and $Z$ be the zero set of $f$. Assume that $Z$ is not empty. 
Then for a compact set $K$ in $U$, there exist positive constants $J$ and $J'$ depending on $f$ and $K$ such that uniformly for all $x\in K$, $|x-Z|^{J}\leq J'|f(x)|$.
It is easy to see that using this, one may construct approximating distinct points $y_1',...,y_k',z_1',...,z_k'\in \mathbb R^n$ (zeroes of a suitable $f$) with the following two properties:
(1) There exist positive constants $J,J'>0$ depending only on $n$ such that
\[
|(y_1,...,y_k,z_1,...,z_k)-(y_1',...,y_k',z_1',...,z_k')| \leq
J\varepsilon^{J'}.
\]
In particular, we have
\[
|y_i-y_i'|\leq J\varepsilon^{J'}
\]
and
\[
|z_i-z_i'|\leq J\varepsilon^{J'}.
\]
(2) $|y_i'-y_j'|=|z_i'-z_j'|$ for every $i,j$.
Thanks to (2), we may choose a Euclidean motion $A:\mathbb R^n\to \mathbb R^n$ so that $A(y^\prime_i)= z^\prime_i$ for each $i$.
Also, thanks to (1), there exists positive constants $J_1, J_2$ depending only on $n$ with
\[
|A(y_i)-A(y_i')|\leq J_1\varepsilon^{J_2}
\]
So it follows that there exist positive constants $J_3, J_4$ depending only on $n$ with
\[
|A(y_i)-z_i|\leq J_3\varepsilon^{J_4}
\]
which is (b). $\Box$

\begin{rmk}
{\rm In \cite{FD3} (see Example~\ref{e:Example1} and Example~\ref{e:Example2} below), we introduced certain slow rotations and translations as $\varepsilon'$
distorted diffeomorphisms from $\mathbb R^n\to \mathbb R^n$. We call them Slow twists and Slides.
The $\varepsilon$-distortions in  Theorem~\ref{t:lemmafive}, Theorem~\ref{t:lemmafiveprime}, Theorem~\ref{t:Theorem2a}, Theorem~\ref{t:Theorem2b}, Theorem~\ref{t:Theorem3},
Theorem~\ref{t:Theorem4} and Theorem~\ref{t:cextensionblock2} are  built using Slow twists and Slides and so Theorem~\ref{t:lemmafive}, Theorem~\ref{t:lemmafiveprime},
Theorem~\ref{t:Theorem2a}, Theorem~\ref{t:Theorem2b}, Theorem~\ref{t:Theorem3}, 
Theorem~\ref{t:Theorem4} and Theorem~\ref{t:cextensionblock2} are also {\bf Alignment of data results} in Euclidean space. We refer the reader to the papers by Werman and his 
collaborators and the many references cited therein for a good insight into the subject of alignment of data in Euclidean space. 
\cite{WW,DSW,PW,AKMSW,WO,BW,OW}.}
\label{r:twistsslides}
\end{rmk}

We are now ready to begin with the proofs of Theorem~\ref{t:Theorem2a}, Theorem~\ref{t:Theorem2b}, Theorem~\ref{t:Theorem3}, 
Theorem~\ref{t:Theorem4} and Theorem~\ref{t:cextensionblock2}. The remainder of the paper will be devoted to developing machinary to this end.

We conclude this section with some additional notation we will use throughout. Henceforth, we will write $c, c', C, C$' ect. to denote positive constants depending only on $D$ and which may denote different constants in different occurrences and we will write $c_K, c_K', C_K, C'_K$ ect. to denote positive constants depending only on $D$ and $K$ and which may denote different constants in different occurrences. We will write $J,J',J_i$ ect. to denote positive 
constants depending only on $n$ and which may denote different constants in different occurrences.
In this paper, all Euclidean motions will be denoted by $A$, $A^*$, $A^{**}, A_i$ ect. which may denote different Euclidean motions in different occurrences. By $A_{\infty}$ we will mean the identity Euclidean motion.
Throughout, $B(x,r)$ will always denote the open ball in $\mathbb R^D$ with center $x$ and radius $r$. Sometimes we will just write $B$ to mean a Ball in $\mathbb R^D$ when we do not need to specify its center and radius. Henceforth, if not stated otherwise, all diffeomorphisms and
 Euclidean motions will be from $\mathbb R^D$ to $\mathbb R^D$.

We begin with a deeper analysis of $\varepsilon'$-distorted diffeomorphisms from $\mathbb R^n\to \mathbb R^n$ . To this end we have:

\section{$\varepsilon'$-distorted diffeomorphisms from $\mathbb R^n\to \mathbb R^n$; further properties and examples: Slow twists and Slides.}
\setcounter{equation}{0}
In this section, we are going to define Slow twists and Slides as refered to first 
in Remark~\ref{r:twistsslides}.

Slow twists are $\varepsilon'$-distorted diffeomorphisms whose argument is a function of distance from the origin. These rotations reduce their speed of rotation for decreasing argument, are non-rigid 
for decreasing argument becoming rigid for increasing argument. See (\ref{e:lemmaone}). Slides are translations which are $\varepsilon'$-distorted diffeomorphisms and satisfy (\ref{e:lemmatwo}).
Theorem~\ref{t:lemmafive}, Theorem~\ref{t:lemmafiveprime}, Theorem~\ref{t:Theorem2a}, Theorem~\ref{t:Theorem2b}, Theorem~\ref{t:Theorem3}, 
Theorem~\ref{t:Theorem4} and Theorem~\ref{t:cextensionblock2} rely on Slow twists and Slides to build them. (See \cite{DSW, NDS} for further applications of this idea). 
We are going to show in the next section how to gradually build $\varepsilon$-distorted diffeomorphisms from $\mathbb R^D\to \mathbb R^D$ from Euclidean motions
$\mathbb R^D\to \mathbb R^D$ using Slow twists and Slides. \footnote{This a natural idea given we recall that
$O(n)$ is by definition the group of isometries of $\mathbb R^n\to\mathbb R^n$ which perserve a fix point ($SO(n)$ is the subgroup of $O(n)$ of orthogonal matrices of determinant 1) and Euclidean motions from 
$\mathbb R^n\to \mathbb R^n$ are the elements of the symmetry group of $\mathbb R^n$, ie all isometries of $\mathbb R^n$.}

\subsection{Three useful properties of $\varepsilon'$-distorted diffeomorphisms from $\mathbb R^n\to \mathbb R^n$.}
We record the following three useful properties of $\varepsilon'$-distorted diffeomorphisms from $\mathbb R^n \to \mathbb R^n$ we will need and which 
follow from (\ref{e:edistortion}).
\begin{itemize}
\item If $\Phi$ is $\varepsilon'$-distorted and $\varepsilon'<\varepsilon''$, then $\Phi$ is $\varepsilon''$-distorted.
\item If $\Phi$ is $\varepsilon'$-distorted, then so is $\Phi^{-1}$.
\item If $\Phi$ and $\Psi$ are $\varepsilon'$-distorted, then the composition map $\Phi o\Psi$ is $J\varepsilon'$-distorted.
\item Suppose $\Phi$ is $\varepsilon'$-distorted. If $\tau$ is a piecewise smooth curve in $\mathbb R^{n}$, then the length of
$\Phi(\tau)$ differs from that of $\tau$ by at most a factor of $(1+\varepsilon')$. Consequently, if $x,x'\in \mathbb R^n$, then $|x-x'|$ and $|\Phi(x)-\Phi(x')|$ differ by at most a factor $(1+\varepsilon)$, ie
we have (\ref{e:edistortion1}):
\[
(1+\varepsilon')^{-1}|x-y|\leq |\Phi(x)-\Phi(y)|\leq |x-y|(1+\varepsilon'),\, x,y\in \mathbb R^n.
\]
\item (\ref{e:edistortion}) together with the fact that $(\Phi'(x))^{T}\Phi'(x)$ is real and symmetric implies that
\[
|(\Phi'(x))^{T}\Phi'(x)-I|\leq J'\varepsilon,\, x\in \mathbb R^n.
\]
This follows from working in an orthonormal basis for which $(\Phi'(x))^{T}\Phi'(x)$ is diagonal.
\end{itemize}

We now meet Slow Twists and Slides.
\subsection{Slow twists.}

\begin{exm}
{\rm Let $\varepsilon'>0$ and $x\in \mathbb R^n$. Let $S(x)$ be the $D\times D$ block-diagonal matrix
\[
\left(
\begin{array}{llllll}
D_1(x) & 0 & 0 & 0 & 0 & 0 \\
0 & D_2(x) & 0 & 0 & 0 & 0 \\
0 & 0 & . & 0 & 0 & 0 \\
0 & 0 & 0 & . & 0 & 0 \\
0 & 0 & 0 & 0 & . & 0 \\
0 & 0 & 0 & 0 & 0 & D_r(x)
\end{array}
\right)
\]
where for each $i$ either $D_i(x)$ is the $1\times 1$ identity matrix or else
\[
D_i(x)=\left(
\begin{array}{ll}
\cos f_i(|x|) & \sin f_i(|x|) \\
-\sin f_i(|x|) & \cos f_i(|x|)
\end{array}
\right)
\]
where $f_i:\mathbb R\to \mathbb R$ are functions satisfying the condition: $t|f_i'(t)|<J_1\varepsilon'$  uniformly for $t\geq 0$. The $1\times 1$ identity matrix is used to compensate for the even/odd size of the matrix. Let 
$\Phi(x)=\Theta^{T}S(\Theta x)$ where $\Theta$
is any fixed matrix in $SO(n)$. Then $\Phi:\mathbb R^n\to \mathbb R^n$ is a $\varepsilon'$-distorted diffeomorphism and we call it a {\it slow twist} (in analogy to rotations). }
\label{e:Example1}
\end{exm}

\subsection{Slides.}
\begin{exm}
{\rm Let $\varepsilon'>0$ and let $g:\mathbb R^{n}\to \mathbb R^{n}$
be a smooth map such that $|g'(t)|<J_2\varepsilon'$ uniformly for $t\in \mathbb R^n$. 
Consider the map $\Phi(t)=t+g(t)$, $t\in \mathbb R^n$. Then $\Phi:\mathbb R^n\to \mathbb R^n$ is a $\varepsilon'$ distorted diffeomorphism. We call the map $\Phi$ a {\it slide} (in analogy to translations). }
\label{e:Example2}
\end{exm}

We now ready to start to build and so to this end we have our next section:

\section{Building $\varepsilon$-distorted diffeomorphisms from Euclidean motions.}
\setcounter{equation}{0}
Our main result of this section is:

\begin{thm}
Let $\varepsilon>0$, $r>0$, $x_1\in \mathbb R^D$, let $B(x_1,r)$ be a ball and let $A$ and $A^*$ be proper Euclidean motions such that
\beq
|A(x_1)-A^*(x_1)|\leq \varepsilon r.
\label{e:emotionsuu}
\eeq
Then there exists a $C\varepsilon$-distorted diffeomorphism $\Phi$ such that $\Phi=A$ in
$B(x_1,\exp(-1/\varepsilon)r)$ and $\Phi=A^*$ outside $B(x_1, r)$.
\label{t:lemmatwists}
\end{thm}

{\bf Proof of Theorem~\ref{t:lemmatwists}}. Assume the hypotheses of the theorem. As a consequence of Slow twists (see Example~\ref{e:Example1}), we have (see also \cite{FD3}), 
that there exists $\eta>0$ depending on $\varepsilon$ for which the following holds. Let $\Theta\in SO(D)$, $r_1,r_2>0$ and let $0<r_1\leq \eta r_2$.
Then, there exists an $\varepsilon$-distorted diffeomorphism $\Phi:\mathbb R^{D}\to \mathbb R^{D}$ such that
\beq
\left\{
\begin{array}{ll}
\Phi(x)=\Theta x, & |x|\leq r_1 \\
\Phi(x)=x, & |x|\geq r_2
\end{array}
\right.
\label{e:lemmaone}
\eeq
Also, as a consequence of Slides (see Example~\ref{e:Example2}), we have (see also \cite{FD3}), 
that there exists $\eta_1>0$ depending on $\varepsilon$ such that the following holds. Let $A:x\to Ex+x_0:\mathbb R^D\to \mathbb R^D$ be a proper Euclidean
motion. Let $r_3,r_4>0$.
Suppose $0<r_3\leq \eta_{1} r_4$ and $|x_0|\leq c\varepsilon r_3$.
Then there exists an $\varepsilon$-distorted diffeomorphism $\Phi_1:\mathbb R^{D}\to \mathbb R^{D}$ such that
\beq
\left \{
\begin{array}{ll}
\Phi_{1}(x)=A(x), & |x|\leq r_3 \\
\Phi_{1}(x)=x, & |x|\geq r_4
\end{array}
\right.
\label{e:lemmatwo}
\eeq

(\ref{e:lemmaone}) and (\ref{e:lemmatwo}) then imply, (see also \cite{FD3}) 
that there exists $\eta_2>0$ depending on $\varepsilon$ such that the following holds. Let $r_5,r_6>0$ with $0<r_5\leq \eta_2 r_6$ and let $x,x'\in \mathbb R^{D}$ with $|x-x'|\leq c\varepsilon r_5$ and 
$|x|\leq r_5$. Then, there exists an $\varepsilon$-distorted diffeomorphism $\Phi:\mathbb R^{D}\to\mathbb R^{D}$ such that 
\beq
\Phi(x)=x'\, {\rm and}\, \Phi(y)=y,\, {\rm for}\, |y|\geq r_6.
\label{e:corollaryone}
\eeq

Theorem~\ref{t:lemmatwists} then follows from (\ref{e:lemmaone}), (\ref{e:lemmatwo}) and (\ref{e:corollaryone}). $\Box$

We are now going to introduce the notation $E$ for a special finite set in $\mathbb R^D$ whose diameter satisfies ${\rm diam}(E)\leq 1$ and its points are well separated (we will define this more precisely in a moment and throughout the paper).
The aim of the next section is to show that for certain such sets $E$, 
we can always construct an improper $\varepsilon$-distorted
diffeomorphism $\Phi:\mathbb R^D\to\mathbb R^D$ such that $\Phi(z)=z$ for each $z\in E$, $\Phi$ agrees with a improper Euclidean motion $A_z:\mathbb R^D\to \mathbb R^D$ in a ball of small enough radius (in particular smaller than the maximum separation distance between points in $E$) and with center $z$ for each $z\in E$ and $\Phi$ agrees with a improper Euclidean motion on all points in $\mathbb R^D$ whose distance to $E$ is large enough. Such sets $E$ we will use to define our special constant $K$ in (~\ref{e:specialk}). 

Our technique will be one of Approximate Reflections. Thus we intoduce our next section of:

\section{Approximate Reflections from $\mathbb R^D$ to $\mathbb R^D$.}
\setcounter{equation}{0}
Suppose that $S$ is a finite subset of a affine hyperplane $H\subset \mathbb R^D$. (So $H$ has dimension $D-1$). Let $A:\mathbb R^D\to \mathbb R^D$ denote reflection through $H$. Then 
$A$ is an improper Euclidean motion and $A(z)=z$ for each $z\in S$. \footnote{For easy understanding: Suppose $D=2$ and $H$ is a line with the set $S$ on the line. Let $A$ denote reflection of the lower half plane to the upper half plane through $S$. Then $A$ is a Euclidean motion and fixes points on $S$ because it is an isometry.} Now suppose that $S$ is again a finite subset of a affine hyperplane $H\subset \mathbb R^D$ and assume that 
we have a map $p:\mathbb R^D \to \mathbb R^D$ with $p(z)$ close to $z$  on $S$. (We will define this more precisely in a moment).  Then we call $p$ an approximate reflection through $H$.
We will now use approximate reflections to construct $\varepsilon$-distorted diffeomorphisms as we have described above.

We proceed and for our set $S$, we will now take a special set $E\subset \mathbb R^D$ as mentioned above which we now define precisely in our main result of this next subsection:
\subsection{Theorem~\ref{t:lemmareflection3}.}
We have:

\begin{thm}
Let $\varepsilon>0$, $0<\tau<1$, $E\subset \mathbb R^D$ be a finite set with ${\rm diam}(E)=1$ and $|z-z'|\geq \tau$ for all $z,z'\in E$ distinct. Assume that $V_D(E)\leq \eta^{D}$ where $0<\eta<c\tau\varepsilon$ for small enough $c$. Here we recall $V_D$ is given by Definition~\ref{d:block}. Then, there exists a $C\varepsilon$-distorted diffeomorphism $\Phi:\mathbb R^D\to \mathbb R^D$ with the following properties:
\begin{itemize}
\item[(a)] $\Phi$ coincides with an improper Euclidean motion on $\left\{x\in \mathbb R^D:\, {\rm dist}(x, E)\geq 20\right\}$.
\item[(b)] $\Phi$ coincides with an improper Euclidean motion $A_z$ on $B(z,\tau/100)$ for each $z\in E$.
\item[(c)] $\Phi(z)=z$ for each $z\in E$.
\end{itemize}
\label{t:lemmareflection3}
\end{thm}

The proof of Theorem~\ref{t:lemmareflection3} relies on two Lemmas, Lemma~\ref{l:lemmareflection1} and Lemma~\ref{l:lemmareflection2}. We begin with:

\begin{lem}
Let $\varepsilon>0$, $0<\tau<1$, $E\subset \mathbb R^D$ be finite set. Assume that ${\rm diam}(E)\leq 1$ and $|z-z'|\geq \tau$ for $z,z'\in E$ distinct. Let $p: \mathbb R^D\to \mathbb R^D$ and assume one 
has $|p(z)-z|\leq \eta$ for all $z\in E$ where 
$\eta<c\varepsilon \tau$ for small enough $c$. \footnote{$p$ is our approximate reflection}. Then, there exists a $\varepsilon$-distorted diffeomorphism
$\Phi$ such that:
\begin{itemize}
\item[(a)] $\Phi(x)=x$ whenever ${\rm dist}(x, E)\geq 10$, $x\in\mathbb R^D$. 
\item[(b)] $\Phi(x)=x+[z-p(z)]$ for $x\in B(z, \tau/10)$, $z\in E$.
\end{itemize}
\label{l:lemmareflection1}
\end{lem}

{\bf Proof:} Let $\theta(y)$ be a smooth cutoff function on $\mathbb R^D$ such that $\theta(y)=1$ for
$|y|\leq 1/10$, $\theta(y)=0$ for $|y|\geq 1/5$ and $|\nabla \theta|\leq C$ on $\mathbb R^D$. 
Let
\[
f(x)=\sum_{z\in E}(z-p(z))\theta\left(\frac{x-z}{\tau}\right),\, x\in \mathbb R^D.
\]
Let $x\in \mathbb R^D$. We observe that if ${\rm dist}(x, E)\geq 10$, then $\frac{|x-z|}{\tau}\geq 1/5$ for each $z\in E$ and so $\theta\left(\frac{x-z}{\tau}\right)=0$
for each $z\in E$. Thus $f(x)=0$ if ${\rm dist}(x, E)\geq 10$. Next if $x\in B(z,\tau/10)$ for each $z\in E$, then for a given $z\in E$, say $z_i$ and for $x\in B(z_i,\tau/10)$, $\frac{|x-z_i|}{\tau}\leq 1/10$ and so $\theta\left(\frac{x-z_i}{\tau}\right)=1$. However since $|z-z'|\geq \tau$ for $z,z'\in E$ distinct, we also have for $x\in B(z_,\tau/10), z\in E,\, z\neq z_i$, $\frac{|x-z|}{\tau}\geq 1/5$ and so $\theta\left(\frac{x-z}{\tau}\right)=0$ for $z\neq z_i$. 
Thus $f(x)=z-p(z)$ for $x\in B(z, \tau/10)$, $z\in E$. Finally 
$|\nabla f|\leq \frac{\eta}{C\tau}<c\varepsilon$ where $c$ is small enough. Then the map $x\to x+f(x),\, x\in \mathbb R^D$ is a slide as per Example~\ref{e:Example2} and thus a $\varepsilon$-distorted diffeomorphism. Thus Lemma~\ref{l:lemmareflection1} holds. $\Box$.

Notice that we do not require $\rho$ to be a Euclidean motion for Lemma~\ref{l:lemmareflection1} to hold. We are now going to show that
under more restrictive conditions on $E$, we can indeed find an improper Euclidean motion $A:\mathbb R^D\to \mathbb R^D$ satisfying
$|A(x)-x|\leq \eta$ for all $x\in E$. This will be Lemma~\ref{l:lemmareflection2} below.

We will need to recall some facts about tensors in $\mathbb R^D$ and wedge and tensor product before we state and prove Lemma~\ref{l:lemmareflection2}.
\subsection{Quick discussion of tensors in $\mathbb R^D$, wedge and tensor product.}
We recall that $\mathbb R^D\otimes \mathbb R^{D'}$ is a subspace of $\mathbb R^{DD'}$ of dimension $DD'$. Here, if $e_1,...e_D$ and $f_1,...f_{D'}$ are the standard basis vectors of $\mathbb R^D$ and $\mathbb R^{D'}$ respectfully, then if $x=(x_1,...x_D)\in \mathbb R^D, \, x_i\in \mathbb R$ and 
$y=(y_1,...,y_D')\in \mathbb R^{D'},\, y_i\in \mathbb R$, $x\otimes y$ is a vector in $\mathbb R^{DD'}$ spanned by the basis vectors $e_i\otimes f_j$ with coefficients $x_iy_j$.\footnote{ For example if $D=2, D'=3$, then $x\otimes y$ is the vector in $\mathbb R^6$ with 6 components $(x_1y_1, x_1y_2, x_1y_3, y_2x_1, y_2x_2, y_2x_3)$ .}
The wedge product of $x,y$, $x\wedge y$ is the antisymmetric tensor product $x\otimes y-y\otimes x$. A tensor in $\mathbb R^D$ consists of components and also basis vectors associated to each component. The number of components of a tensor in $\mathbb R^D$ need not be $D$.
The rank of a tensor in $\mathbb R^D$ is the minimum number of basis vectors in $\mathbb R^D$ associated to each component of the tensor. ($\mathbb R^D$ realized as vectors $(p_1,...,p_D)$ where $p_i\in \mathbb R$ are rank-1 tensors given each component $p_i$ has 1 basis vector $e_i$ associated to it.) A rank-l tensor in $\mathbb R^D$ has associated to each of its components $l$ basis vectors out of $e_1,e_2,e_3,...e_D$. Real numbers are 0-rank tensors.
Let now $v_1, ...v_l\in \mathbb R^D$. Writing
$v_1=a_{11}e_1+...+a_{l1}e_l,v_2=a_{12}e_1+...+a_{l2}, ..,v_l=a_{1l}e_1+...+a_{ll}e_l$ the following holds:
\begin{itemize}
\item $v_1\wedge...\wedge v_l=({\rm det}A)e_1\wedge e_2...\wedge e_l$ where $A$ is the matrix
\[
\left(
\begin{array}{llll}
a_{11} & a_{12} &...& a_{1l} \\
a_{21} & a_{22} &...& a_{2l} \\
.& .& . & . \\
a_{l1}& a_{l2} & ... & a_{ll}
\end{array}
\right).
\]
\item $|v_1\wedge v_2...\wedge v_l|^2={\rm det}(A')={\rm vol}_l(v_1,v_2...v_l)$.
Here, $A'$ is the matrix
\[
\left(
\begin{array}{llll}
v_1.v_1 & v_1.v_2 &...& v_1.v_l \\
v_2.v_1& v_2.v_2 &...& v_2.v_l \\
.& .& . & . \\
v_l.v_1& v_l.v_2 & ... & v_l.v_l
\end{array}
\right)
\]
and \[
{\rm vol}_l(v_1,v_2...v_l)=\left\{\sum_{i=1}^{l}c_iv_i:\, 0\leq c_i\leq 1,\, c_i\in \mathbb R\right\}
\]
is the l-volume of the parallelepiped determined by $v_1,...v_l$. $|.|$ is here understod as the rotationally invariant norm on alternating tensors of any rank.
\end{itemize}

We are ready for:
\subsection{ Lemma~\ref{l:lemmareflection2}.}
Given a finite set $S\subset \mathbb R^D$, we recall the definition of $V_D(S)$. See Definition~\ref{d:block}. We are now ready for:

\begin{lem}
Let $0<\eta<1$ and $E\subset \mathbb R^D$ finite with ${\rm diam}(E)=1$. Assume that $V_{D}(E)\leq \eta^D$. Then, there exists an improper Euclidean motion $A:\mathbb R^D\to \mathbb R^D$ such that
\beq
|A(z)-z|\leq C\eta, \forall z\in E.
\eeq
\label{l:lemmareflection2}
\end{lem}

{\bf Proof:}\, We have $V_1(E)=1$ and $V_D(E)\leq \eta^D$. Hence, there exists $l$ with $2\leq l\leq D$ such that
$V_{l-1}(E)>\eta^{l-1}$ but $V_{l}(E)\leq \eta^l$. Fix such a $l$. Then there exists a $(l-1)$ simplex with vertices
$z_0,...,z_{l-1}\in E$ and with $(l-1)$ dimensional volume $>\eta^{l-1}$. Fix $z_0,...,z_{l-1}$. Without loss of generality, we may suppose $z_0=0$. Then
\[
|z_1\wedge,...,\wedge z_{l-1}|>c\eta^{l-1}
\]
yet
\[
|z_1\wedge,...,\wedge z_{l-1}\wedge z|\leq C\eta^{l}
\]
for any $z\in E$. Now,
\[
|z_1\wedge,...,\wedge z_{l-1}\wedge z|=|\pi z||z_1\wedge...\wedge z_{l-1}|
\]
where $\pi$ denotes the orthogonal projection from $\mathbb R^D$ onto the space of vectors orthogonal to $z_1,...,z_{l-1}$. Consequently, we have for $z\in E$,
\[
C\eta^l\geq |z_1\wedge,...,\wedge z_{l-1}\wedge z|=|\pi z||z_1\wedge,...,\wedge z_{l-1}|\geq c\eta^{l-1}|\pi z|.
\]

We deduce that we have $|\pi z|\leq C\eta$ for any $z\in E$. Equivalently, we have shown that every $z\in E$ lies within a distance $C\eta$ from ${\rm span}\left\{z_1,...,z_{l-1}\right\}$.
This span has dimension $l-1\leq D-1$. Letting $H$ be the hyperplane containing that span and letting $A$ denote the reflection through $H$, we see that ${\rm dist}(z,H)\leq C\eta$. Hence,
\[
|A(z)-z|\leq C\eta,\, \forall z\in E.
\]
Since $A$ is an improper Euclidean motion, the proof of the Lemma is complete. $\Box$.

We are now ready to give the proof of Theorem~\ref{t:lemmareflection3}.
\medskip

{\bf Proof:}\, By Lemma~\ref{l:lemmareflection2}, there exists an improper Euclidean motion $A:\mathbb R^D\to \mathbb R^D$ such that for $z\in E$,
\[
|A(z)-z|\leq C\eta.
\]
Hence, Lemma~\ref{l:lemmareflection1} applies with $C\eta$ in place of $\eta$. Let $\Psi$ be an $\varepsilon$-distorted diffeomorphism as in the conclusion of Lemma~\ref{l:lemmareflection1}. 
Then it is easy to check that $\Phi:=\Psi o A $ is a $C\varepsilon$-distorted diffeomorphism and satisfies the conclusion we desire. We are done. 
$\Box$
\medskip

\section{Approximation by Euclidean Motions.}
\setcounter{equation}{0}
We now want an analogy of Theorem~\ref{t:lemmareflection3} for a proper $\varepsilon$- diffeomorphisms (see Theorem~\ref{t:lemmareflectionex4} below) and it is here that we meet the constant $K$. See (\ref{e:specialk}) 
below.
In order to do this, we need more machinery. To begin, it will be necessary first
to study pointwise approximation of $\varepsilon$ diffeomorphisms by given Euclidean motions. It follows from \cite{DF4} 
that given $\varepsilon''>0$ (small enough and depending only on $D$) and $\Phi:\mathbb R^D\to \mathbb R^D$ a $\varepsilon''$-distorted diffeomorphism, there exists a Euclidean motion $A:\mathbb R^D\to \mathbb R^D$ with 
$|\Phi(x)-A(x)|\leq C\varepsilon''$ for $x\in B(0,10)$. 
Actually, using the well-known John-Nirenberg inequality, in \cite{DF4} we proved a lot more, namely a BMO theorem for $\varepsilon''$-distorted diffeomorphisms which is in the next subsection. 

\subsection{BMO theorem for $\varepsilon'$-distorted diffeomorphisms.} 
{\bf BMO theorem for $\varepsilon''$-distorted diffeomorphisms}:\, Let $\varepsilon''>0$ be a small enough positive number depending only on $D$,  $\Phi:\mathbb R^D\to \mathbb R^D$  a $\varepsilon''$-distorted diffeomorphism and let $B\in \mathbb R^D$ be a ball. There exists $T=T_B\in O(D)$ and 
$C>0$ such that for all $\lambda\geq 1$, we have
\beq
{\rm vol}\left\{x\in B:|\Phi'(x)-T(x)|>C\varepsilon''\lambda\right\}\leq \exp(-\lambda){\rm vol}(B)
\eeq
{\rm and slow twists in Example~\ref{e:Example1} show that the estimate above is sharp. The set 
\beq
\left\{x\in B:|\Phi'(x)-T(x)|> C\varepsilon''\lambda\right\}
\eeq may well be small in a more refined sense but we did not pursue this investigation in \cite{DF4}}.

\subsection{Approximation by Euclidean Motions; Proper and Improper.}
From Definition~\ref{d:block}, we recall that $V_l(z_0,...,z_l)$ denotes the $l$-dimensional volume of the $l$-simplex with vertices at $z_0,...,z_l$. We now have:

\begin{thm}
\item[(a)] Let $\varepsilon>0$. Let $\Phi:\mathbb R^D\to \mathbb R^D$ be a $\varepsilon$-distorted diffeomorphism. Let $z\in \mathbb R^D$ and $r>0$ be given.  Then, there exists an Euclidean motion $A=A_B$ with  
$B=B(z,r)$,  such that for $x\in B(z,r)$,
\begin{itemize}
\item[(1)] $|\Phi(x)-A(x)|\leq C\varepsilon r$.
\item[(2)] Moreover, $A$ is proper iff $\Phi$ is proper.
\end{itemize}
\item[(b)] Let $x_0,...,x_D\in \mathbb R^D$ with ${\rm diam}\left\{x_0,...,x_D\right\}\leq 1$ and $V_D(x_0,...,x_D)\geq \eta^D$ where
$0<\eta<1$ and let $0<\delta<c'\eta^D$ for a small enough $c'$. Let $\Phi:\mathbb R^D\to\mathbb R^D$ be a $\delta$-distorted diffeomorphism. Finally let $T$ be the one and only one affine map that agrees with
$\Phi$ on $\left\{x_0,...,x_D\right\}$. (We recall that the existence and uniqueness of such $T$ follows from \cite{ATV}, T may not be invertible). Then $\Phi$ is proper iff $T$ is proper.
\label{t:theorememotionapprox}
\end{thm}

{\bf Proof:}\, We begin with part (a).  Without loss of generality, we may assume that $B(z,r)=B(0,1)$ and $\Phi(0)=0$. Let $e_1,...,e_D\in \mathbb R^D$ be unit vectors. Then, $|\Phi(e_i)|=|\Phi(e_i)-\Phi(0)|$. Hence, for each $i$,
\[
(1+\varepsilon)^{-1}\leq |\Phi(e_i)|\leq (1+\varepsilon).
\]
Also for $i\neq j$,
\[
(1+\varepsilon)^{-1}\sqrt{2}\leq |\Phi(e_i)-\Phi(0)|\leq (1+\varepsilon)\sqrt{2}.
\]
Hence,
\[
\Phi(e_i).\Phi(e_j) =1/2\left(|\Phi(e_i)|^2+\Phi(e_j)|^2-|\Phi(e_i)-\Phi(e_i)|^2\right)
\]
satisfies
\[
|\Phi(e_i).\Phi(e_j)-\delta_{ij}|\leq C\varepsilon
\]
for all $i,j$ where $\delta_{ij}$ denotes the Kroncker delta and "." denotes the Euclidean dot product.
Applying the Gram-Schmidt process to $\Phi(e_1),....,\Phi(e_D)$, we obtain orthonormal vectors $e_1^*,...,e_D^*\in \mathbb R^D$ such that $|\Phi (e_i)-e_i^*|\leq C\varepsilon$ for each $i$. Using Theorem~\ref{t:Theorem 1}, we let $A$ be the (proper or improper) rotation such that $Ae_i=e_i^*$ for each $i$. Then $\Phi^{**}:=A^{-1}o\Phi$ is an $\varepsilon$-distorted diffeomorphism,
$\Phi^{**}(0)=0$ and $|\Phi^{**}(e_i)-e_i|\leq C\varepsilon$ for each $i$. Now let $x=(x_1,...,x_D)\in B(0,1)$ and let $y=(y_1,...,y_D)=\Phi^{**}(x)$. Then $2x_i=1+|x-0|^2-|x-e_i|^2$ and also
$2y_i=1+|y-0|^2-|y-e_i|^2$ for each i. Hence, by the above-noted properties of $\Phi^{**}$, we have $|y_i-x_i|\leq C\varepsilon$. Then, $|\Phi^{**}(x)-x|\leq C\varepsilon$ for all $x\in B(0,1)$, ie,
$|\Phi(x)-A(x)|\leq C\varepsilon$ for all $x\in B(0,1)$. Thus, we have proved (1) but not yet (2). For each $(z,r)$, (1) provides an Euclidean motion $A_{(z,r)}$ such that $|\Phi(x)-A_{(z,r)}(x)|\leq C\varepsilon r$ for $x\in B(z,r)$.
Now for $r$ small enough, we have using the mean value theorem for vector valued functions and the substitution rule with Jacobian determinants as expansions of volumes,
\[
|\Phi(x)-[\Phi(z)+\Phi'(z)(x-z)]|\leq C\varepsilon r, \, x\in B(z,r).
\]
Hence,
\[
|A_{(z,r)}(x)-[\Phi(z)+\Phi'(z)(x-z)]|\leq C\varepsilon r, \, x\in B(z,r).
\]
Thus we have established for small $r$ that $A_{(z,r)}$ is proper iff ${\rm det}\Phi'(z)>0$ ie, iff $\Phi$ is proper.
Observe that $|\Phi(x)-A_{(z,r)}(x)|\leq C\varepsilon r$ for $x\in B(z,r)$ iff $|\Phi(x)-A_{(z,r/2)}(x)|\leq C\varepsilon r$ for $x\in B(z,r/2)$. Thus
$|A_{(z,r)}-A_{(z,r/2)}(x)|\leq C\varepsilon r$ for $x\in B(z,r/2)$.
Hence $A_{(z,r)}$ is proper iff $A_{(z,r/2)}$ is proper. Thus we may deduce that for all $r$, $A_{(z,r)}$ is proper iff $\Phi$ is proper. This completes the proof of (2) and part (a) of Theorem~\ref{t:theorememotionapprox}.
\medskip

We now prove part (b). Without loss of generality, we may assume that $x=0$ and $\Phi(x_0)=0$. Then $T$ is linear, not just affine. By Theorem~\ref{t:theorememotionapprox}, there exists a
Euclidean motion $A_{(0,1)}$ such that
\[
|\Phi(x)-A_{(0,1)}(x))|\leq C\delta
\]
for all $x\in B(0,1)$ and $\Phi$ is proper iff $A_{(0,1)}$ is proper. We know that
\[
|Tx_i-A_{(0,1)}(x_i)|\leq C\delta, i=0,1,...D
\]
since $Tx_i=\Phi(x_i)$ and also since $x_i\in B(0,1)=B(x_0,1)$. (The later follows because ${\rm diam}\left\{x_0,...,x_D\right\}\leq 1$). In particular, $|A_{(0,1)}(0)|\leq C\delta$ since $x_0=0$.
Hence,
\[
|Tx_i-[A_{(0,1)}(x_i)-A_{(0,1)}(0)]|\leq C'\delta
\]
for $i=1,...,D$.
Now, the map $x\mapsto Ax:=A_{(0,1)}(x)-A_{(0,1)}(0)$ is a proper or improper rotation and ${\rm det}(A)>0$ iff
$A_{(0,1)}$ is proper iff $\Phi$ is proper. Thus we have the following:
\begin{itemize}
\item $|(T-A)x_i|\leq C'\delta$, $i=1,...,D$.
\item $|x_1\wedge...\wedge x_D|\geq c\eta^D$ (from the discussion in Section 6.2).
\item ${\rm det}A>0$ iff $\Phi$ is proper.
\item $A$ is a proper or improper rotation.
\end{itemize}

Now let $L$ be the linear map that sends the $i$th unit vector $e_i$ to $x_i$. Then the entries of $L$ are at most $1$ in absolute value since each $x_i$ belongs to $B(0,1)$. Letting $|.|$ be understod as the rotationally invariant norm on alternating tensors of any rank or the operator norm on the space of $D\times D$ square matrices viewed as linear operators from $\mathbb R^D \to \mathbb R^D$ induced from the Euclidean norm on $\mathbb R^D$ or
otherwise the Euclidean norm on $\mathbb R^D$, we have from the discussion in Section 6.2 and the above,
\[
|{\rm det}(L)|=|x_1\wedge...\wedge x_D|\geq c\eta^D.
\]
Hence by Cramers rule, $|L^{-1}|=\left|\frac{{\rm Adjoint}(L)}{{\rm det}(L)}\right|\geq C\eta^{-D}$. We now have for each $i$,
\[
|(T-A)Le_i|=|(T-A)x_i|\leq C'\delta.
\]
Hence,
\[
|(T-A)L|\leq C''\delta
\]
and thus (recall $(T-A)L$ and $L^{-1}$ are square),
\[
|T-A|=|(T-A)LL^{-1}|\leq |(T-A)L||L^{-1}|\leq C\delta \eta^{-D}.
\]
Since $A$ is a (proper or improper) rotation, it follows that (recall $(T-A)L$ and $L^{-1}$ are square),
\[
|TA^{-1}-I|=|(T-A)A^{-1}|\leq |(T-A)||A^{-1}|\leq C\delta\eta^{-D}.\]

Therefore if $\delta\eta^{-D}\leq c'$ for small enough $c'$, then $TA^{-1}$ lies in a
small neighborhood of $I$ and therefore ${\rm det}(TA^{-1})>0$. Hence ${\rm det}T$ and ${\rm det}(A)$ have the same sign. Thus, ${\rm det}(T)>0$ iff $\Phi$ is proper. So we have proved (b) and the theorem.$\Box$.

\subsection{Proof of Theorem~\ref{t:cextensionblock2}.}
We have:

{\bf Proof of Theorem~\ref{t:cextensionblock2}}:

{\bf Proof}:\, The theorem follows from Theorem~\ref{t:theorememotionapprox} part (b) using Definition~\ref{d:block} with a rescalling
${\rm diam}\left\{x_0,...,x_D\right\}\leq 1$ and $V_D(x_0,...,x_D)\geq \eta^{D}{\rm diam}\left\{x_0,...,x_D\right\}. $ $\Box$.

We note that the same argument in the  proof of Theorem~\ref{t:cextensionblock2} yields:

\begin{cor}
Let $\phi:E\to \mathbb R^D$ where $E\subset \mathbb R^D$ is finite. Assume that $\phi$ has a positive (resp. negative) $0<\eta<1$ block and let $0<\delta<c\eta^D$ for small enough $c$. 
Suppose $\phi$ extends to a $\delta$ distorted diffeomorphism $\Phi$.
Then, $\Phi$ is proper (resp. improper).
\label{c:cextensionblock1}
\end{cor}

\section{The constant $K$ and Theorem~\ref{t:lemmareflection3} for proper $\varepsilon$-distorted diffeomorphisms.}
\setcounter{equation}{0}
We are now able to prove an analogy of Theorem~\ref{t:lemmareflection3} for proper $\varepsilon$-distorted diffeomorphisms, namely Theorem~\ref{t:lemmareflectionex4}.
Two important ingredients we will need to do this will be to assume a separation condition on the points of $E$ as in Theorem~\ref{t:lemmareflection3} and also a condition on the cardinality of $E$ which will be our constant $K$, see (\ref{e:specialk}) below.

Thus we have as our main result in this section:

\begin{thm}
Let $\phi:E\to \mathbb R^D$ with $E\subset \mathbb R^D$ finite. Let $\varepsilon>0$ and $0<\tau<1$. 
We make the following assumptions:
\begin{itemize}
\item Assumption on parameters:
\begin{itemize}
\item Let $0<\eta<c\varepsilon \tau$ for small enough $c$.
\item Let $C_K\delta^{1/\rho_k}\tau^{-1}<{\rm min}(\varepsilon,\eta^D)$ for some large enough $C_K>0$ and $\rho_K>0$, the later also depending only on 
$D$ and $K$.
\end{itemize}
\item Assumptions on E: ${\rm diam}(E)=1$,\, $|x-y|\geq \tau$, for any $x,y\in E$ distinct, 
\beq {\rm card}(E)\leq K.
\label{e:specialk}
\eeq
\item Assumption on $\phi$: $\phi$ has no negative $\eta$-blocks and
\[
(1+\delta)^{-1}|x-y|\leq |\phi(x)-\phi(y)|\leq (1+\delta)|x-y|,\, x, y\in E.
\]
\end{itemize}
Then, there exists a proper $C\varepsilon$-distorted diffeomorphism $\Phi:\mathbb R^D\to \mathbb R^D$ with the following
properties:
\begin{itemize}
\item $\Phi=\phi$ on $E$.
\item $\Phi$ agrees with an Euclidean motion $A_{\infty}$ on
$\left\{x\in \mathbb R^D:\, {\rm dist}(x,E)\geq 1000\right\}.$
\item For each $z\in E$, $\Phi$ agrees with a Euclidean motion $A_z$ on $B(z,\tau/1000)$.
\label{t:lemmareflectionex4}
\end{itemize}
\end{thm}

We remark that it follows immediately from the theorem that $A_z=\Phi$ for each $z\in E$ (if $z\in E$, then trivially $z\in B(z,\tau/1000)$ and also $\Phi=\phi$ for each $z\in E$ and so $\Phi=\phi=A_z$ on $E$. 

\subsection{Auxillary Lemmas.}
Theorem~\ref{t:lemmareflectionex4} will follow from 3 lemmas below namely Lemma~\ref{l:lemmareflectionex1}, Lemma~\ref{l:lemmareflectionex2} and Lemma~\ref{l:lemmareflectionex3}.

We begin with:
\begin{lem}
Let $\phi:E\to \mathbb R^D$ with $E\subset \mathbb R^D$ finite and let $0<\tau<1$.
We make the following assumptions:
\begin{itemize}
\item Assumption on E: ${\rm diam}(E)\leq 1$,\, $|x-y|\geq \tau$, for any $x,y\in E$ distinct, ${\rm card}(E)\leq K$.
\item Assumption on $\delta$: $0<\delta\leq c_K\tau^{\rho_{K}}$ for some $c_K>0$ (small enough) and $\rho_K>0$
(large enough) the later depending only on $K,D$.
\item Assumption on $\phi$:
\[
(1+\delta)^{-1}|x-y|\leq |\phi(x)-\phi(y)|\leq (1+\delta)|x-y|,\, x, y\in E.
\]
\end{itemize}

Then, there exists a $C_K\delta^{1/\rho_K}\tau^{-1}$ distorted diffeomorphism $\Phi:\mathbb R^D\to \mathbb R^D$ with the following
properties:
\begin{itemize}
\item $\Phi=\phi$ on $E$.
\item $\Phi$ agrees with an Euclidean motion $A_{\infty}$ on
$\left\{x\in \mathbb R^D:\, {\rm dist}(x,E)\geq 100\right\}.$
\item For each $z\in E$, $\Phi$ agrees with a Euclidean motion $A_z$ on $B(z,\tau/100)$.
\label{l:lemmareflectionex1}
\end{itemize}
Note that we did not say whether $\Phi $ is proper or improper.
\end{lem}

{\bf Proof:} Using the Lojasiewicz inequality, (see Section 3.6), there exists a Euclidean motion $A$ for which we have
\[
|\phi(x)-A(x)|\leq C_K\delta^{1/\rho_K},\, x\in E.
\]
Without loss of generality, we may replace $\phi$ by $\phi^*:=\phi o A^{-1}$. Hence, we may suppose that
\[
|\phi(x)-x|\leq C_k\delta^{1/\rho_{k}}, \, x\in E.
\]
Now we will employ a similar technique to the proof of Lemma~\ref{l:lemmareflection1}.

Let $\theta(y)$ be a smooth cut off function on $\mathbb R^D$ such that $\theta(y)=1$ for $|y|\leq 1/100$, $\theta(y)=0$
for $|y|\geq 1/50$ and with $|\nabla\theta(y)|\leq C$ for all $y$. Then set
\[
f(x)=\sum_{z\in E}(\phi(z)-z)\theta(x-z/\tau),\, x\in \mathbb R^D.
\]
The summands are smooth and have pairwise disjoint supports and thus $f$ is smooth. As in the proof of Lemma~\ref{l:lemmareflection1}, 
$f(x)=0$ for ${\rm dist}(x,E)\geq 100$, $f(x)=\phi(z)-z$ for $x\in B(z,\frac{\tau}{100})$, $z\in E$ and $|\nabla f(x)|\leq C_k\delta^{1/\rho_k}C\tau^{-1}$.
If $CC_k\delta^{1/\rho_k}\tau^{-1}$ is small enough the map $\Phi(x)=f(x)+x$ is a slide and thus $\Phi$ is a $C_k\delta^{1/\rho_k}\tau^{-1}$ distorted diffeomorphism and has all the desired properties. Thus, we are done. $\Box$.

We now worry about whether the map $\Phi$ in Lemma~\ref{l:lemmareflectionex1} is proper or improper. Thus we have:

\begin{lem}
Suppose $E\subset \mathbb R^D$ is finite. Let $\varepsilon>0$ and $0<\tau<1$. Suppose that ${\rm diam}(E)=1$, $|x-y|\geq \tau$ for $x,y\in E$ distinct and ${\rm card}(E)\leq K$. 
Let $0<\eta<c\varepsilon \tau$ for small enough $c$ and let $C_k\delta^{1/\rho_k}\leq \varepsilon \tau$ for large enough $C_K$ and $\rho_K$ depending on $D$ and $K$. Suppose $V_D(E)\leq \eta^D$ where 
$V_D$ is as in Definition~\ref{d:block}. Let $\phi:E\to \mathbb R^D$ and suppose that
\[
(1+\delta)^{-1}|x-y|\leq |\phi(x)-\phi(y)|\leq (1+\delta)|x-y|,\, x,y\in E.
\]
Then, there exists a proper $C\varepsilon$-distorted
diffeomorphism $\Phi:\mathbb R^D\to \mathbb R^D$ with the following properties:
\begin{itemize}
\item $\Phi=\phi$ on $E$.
\item $\Phi$ agrees with a proper Euclidean motion $A_{\infty}$ on
$\left\{x\in \mathbb R^D:\, {\rm dist}(x,E)\geq 1000\right\}.$
\item For each $z\in E$, $\Phi$ agrees with a proper Euclidean motion $A_z$ on $B(z, \tau/1000)$.
\end{itemize}
\label{l:lemmareflectionex2}
\end{lem}

{\bf Proof:}\, Start with $\Phi$ from Lemma~\ref{l:lemmareflectionex1}. If $\Phi$ is proper, then we are done. (Note that
$C_K\delta^{1/\rho_K}\tau^{-1}<\varepsilon$.). If $\Phi$ is improper, then Theorem~\ref{t:lemmareflection3} applies; letting
$\Psi$ be as in Theorem~\ref{t:lemmareflection3}, we see that $\Phi o\Psi$ satisfies all the assertions of Lemma~\ref{l:lemmareflectionex2}. $\Box$.

Finally we have our third lemma:
\begin{lem} Let $E\in \mathbb R^D$, with ${\rm diam}(E)=1$. Let $0<\tau<1$ and suppose $|x-y|\geq \tau$ for $x,y\in E$ distinct. Let $0<\eta<1$ and suppose 
$V_D(E)\geq \eta^D$ and ${\rm card}(E)\leq K$. Let $C_K\delta^{1/\rho_k}\tau^{-1}<{\rm min}(\varepsilon, \eta^D)<1$ for large enough $C_K, \rho_K$ depending on
$D$ and $K$. Let $\phi:E\to\mathbb R^D$ and assume that
\[
(1+\delta)^{-1}|x-y|\leq |\phi(x)-\phi(y)|\leq (1+\delta)|x-y|
\]
for $x,y\in E$ and that $\phi$ has no negative
$\eta$ blocks.
Then, there exists a proper $C\varepsilon$ diffeomorphism $\Phi$ with the following properties.
\begin{itemize}
\item $\Phi=\phi$ on $E$.
\item $\Phi$ agrees with a proper Euclidean motion $A_{\infty}$ on
\[
\left\{x\in \mathbb R^D:\, {\rm dist}(x, E)\geq 1000\right\}.
\]
\item For each $z\in E$, $\Phi$ agrees with a proper Euclidean motion $A_z$ on $B(z, \tau/1000)$.
\end{itemize}
\label{l:lemmareflectionex3}
\end{lem}

{\bf Proof:}\, We apply Lemma~\ref{l:lemmareflectionex1}. The map $\Phi$ in Lemma~\ref{l:lemmareflectionex1} is a $C_K\delta^{1/\rho_K}\tau^{-1}$ distorted diffeomorphism; hence is a $C\varepsilon$-distorted diffeomorphism.
If $\Phi$ is proper, then it satisfies all the conditions needed
and we are done. Thus let us check $\Phi$ is proper. By hypothesis, we can find $z_1,...,z_D\in E$ such that
\[
V_D(z_0,...,z_D)\geq \eta^D.
\]
Let $T$ be the one and only affine map that agrees with $\phi$ on $\left\{z_0,...,z_D\right\}$. See \cite{ATV}. Since $\phi$ has no negative $\eta$ blocks (by hypothesis), we know that $T$ is proper. Applying Theorem~\ref{t:theorememotionapprox} (b) with $\delta$ replaced by
$C_K\delta^{1/\rho_K}\tau^{-1}$, we find that $\Phi $ is proper as needed. Note that Theorem~\ref{t:theorememotionapprox}
applies here since we assumed that $C_K\delta^{1/\rho_K}\tau^{-1}<\eta^D$
for large enough $C_K$ and $\rho_K$ depending only on $K$ and $D$. The proof of Lemma~\ref{l:lemmareflectionex3} is complete. $\Box$.
\medskip

Combining Lemma~\ref{l:lemmareflectionex2} and Lemma~\ref{l:lemmareflectionex3} we are able to give the proof of
Theorem~\ref{t:lemmareflectionex4}.
\medskip

{\bf Proof}\, If $V_D(E)\leq \eta^D$, then Theorem~\ref{t:lemmareflectionex4} follows from Lemma~\ref{l:lemmareflectionex2}. If instead, $V_D(E)>\eta^D$, then Theorem~\ref{t:lemmareflectionex4} follows from Lemma~\ref{l:lemmareflectionex3} $\Box$.

\section{Lets Glue: The Gluing Theorem.}
\setcounter{equation}{0}

Given a finite $E$ with some special geometry and a $\delta$ distortion $\phi$ on $E$, we have investigated in detail up to now how to produce smooth $\varepsilon$-distortions which agree with $\phi$ on the set $E$ and which agree 
with Euclidean motions inside and outside different sets in $\mathbb R^D$. We need now to "Glue" these results together. This is the subject of this section.
\subsection{The Gluing Theorem.}
We prove:

\begin{thm}
Let $E$ be finite, $\varepsilon>0$, $0<\tau<1$, $\phi:E\to \mathbb R^D$ and suppose $|x-y|\geq \tau>0$ for $x,y\in E$ distinct. Suppose also that
\[
1/2|x-y|\leq |\phi(x)-\phi(y)|\leq 2|x-y|
\]
for $x,y\in E$ distinct. For $i=1,...,4$ and $z\in E$, define
\[
B_i(z)=B\left(z,\exp\left((i-5)/\varepsilon\right)\tau\right).
\]
For each $z\in E$, suppose we are given a $C\varepsilon$-distorted diffeomorphism $\Phi_z$ such that $\Phi_z(z)=\phi(z)$ on $E$ and $\Phi_z$ agrees with a proper Euclidean motion $A_z$ outside $B_{1}(z)$ for each $z\in E$.
Moreover, suppose we are given a $C\varepsilon$-distorted diffeomorphism
$\Psi$ such that $\phi=\Psi$ on $E$ and $\Psi$ agrees with a proper Euclidean motion $A_{z}^*$ in $B_4(z)$ for each $z\in E$. Then there exists a $C’\varepsilon$-distorted diffeomorphism $\Phi$ such that:
\begin{itemize}
\item $\Phi=\Phi_z$ in $B_2(z)$ for $z\in E$ (in particular $\Phi=\phi$ on $E$) and
\item $\Phi=\Psi$ outside $\cup_{z\in E}B_3(z)$.
\end{itemize}
\label{t:lemmagl}
\end{thm}

{\bf Proof:}\, We first investigate how well $A_z(z)$ approximates $A_z^*(z)$. Let $z\in E$. Then $A_z^*(z)=\Psi(z)=\phi(z)$ since $z\in B_4(z)$. Moreover, for any $x\in \mathbb R^D$ such that 
$|x-z|=\exp(-4/z)\tau$, we have $x\notin B_1(z)$, hence
$\Phi _z(x)=A_z(x). $ We recall that $\Phi_z$ is a $C\varepsilon$ diffeomorphism and that $\Phi_z(z)=\phi(z)$. Thus,
\[
(1+C\varepsilon)^{-1}|x-z|\leq |\Phi_z(x)-\Phi_z(z)|\leq (1+C\varepsilon)|x-z|
\]
ie,
\[
(1+C\varepsilon)^{-1}\exp(-4/\varepsilon)\tau\leq |A_z(x)-\phi(z)|\leq (1+C\varepsilon)\exp(-4/\varepsilon)\tau. \]
This holds whenever $|x-z|=\exp(-4/\varepsilon)\tau.$ Since $A_z$ is an Euclidean motion, it follows that
\[
|A_z(z)-\phi(z)|\leq C\varepsilon \exp(-4/\varepsilon)\tau.
\]
Recalling that $A_z^*=\phi(z)$, we conclude that for $z\in E$,
\[
|A_z(z)-A_z^*(z)|\leq C \varepsilon \exp(-4/\varepsilon)\tau.
\]
Also, both $A_z$ and $A_z^*$ are proper Euclidean motions. Applying, Theorem~\ref{t:lemmatwists}, we obtain for each $z\in E$, a $C\varepsilon$ diffeomorphism $\Phi_z^*$ such that:
\begin{itemize}
\item $\Phi_z^*$ agrees with $A_z$ on $B_2(z)$.
\item $\Phi_z^*$ agrees with $A_z^*$ outside $B_3(z)$.
\end{itemize}
Let us define a map $\Phi:\mathbb R^D \to \mathbb R^D$ in overlapping regions as follows:
\begin{itemize}
\item $\Phi=\Phi_z$ in $B_2(z)$ for $z\in E$.
\item $\Phi=\Phi_z^*$ in $B_4(z)\setminus B_1(z),\, z\in E$.
\item $\Phi=\Psi$ in $\mathbb R^D\setminus \cup_{z\in E} B_3(z)$.
\end{itemize}
Let us check that the above definitions of $\Phi$ in overlapping regions are mutually consistent.
\begin{itemize}
\item On $B_2(z)\cap [B_4(z')\setminus B_1(z')], z,z'\in E$: To have a non empty intersection, we must have $z'=z$ (since otherwise $|z-z'|\geq \tau$). In the region in question,
$\Phi_z^*=A_z$ (since we are in $B_2(z)$)=$\Phi_z$ (since we are outside $B_1(z)$).
\item On $[B_4 (z)\setminus B_1(z)]\cap [\mathbb R^D\setminus \cup_{z'\in E}B_3(z')],\, z\in E.$ $\Psi=A_z^*$(since we are in $B_4(z))=\Phi_z^*$(since we are outside $B_3(z)$).
\item Note that the balls $B_2(z),\, z\in E$ are pairwise disjoint as are the regions $B_4(z)\setminus B_1(z)$, $z\in E$ since $|z-z'|\geq \tau$ for $z,z'\in E$ distinct.
\end{itemize}
Moreover,
$B_2(z)\cap [\mathbb R^D\setminus \cup_{z'\in E} B_3(z')]=\emptyset.$ Thus, we have already discussed all the non empty intersections of the various regions in which $\Phi$ was defined. This completes the verification that $\Phi$ is defined consistently.

Since $\Psi$, $\Phi_z$, $\Phi_z^*$ (each $z\in E$) are $C\varepsilon$-distorted diffeomorphisms, we conclude that
$\Phi:\mathbb R^D\to \mathbb R^D$ is a smooth map and that
\[
(1+C'\varepsilon)^{-1}\leq (\Phi'(x)^T(\Phi'(x))\leq 1+C'\varepsilon,\, x\in \mathbb R^D.
\]
We have also $\Phi=\Phi_z$ on $B_2(z)$ for each $z\in E$ and $\Phi=\Psi$ outside $\cup_{z\in E}B_3(z)$ by definition of $\Phi$.
To complete the proof of the Gluing theorem, it remains only to check that $\Phi:\mathbb R^D\to \mathbb R^D$ is one to one and onto. To see this, we argue as follows. Recall that the $A_z$ and $A_z^*$ are Euclidean motions and that
\[
|A_z-A_z^*|\leq C\varepsilon\exp(-4/\varepsilon)\tau=C\varepsilon{\rm radius}(B_1(z)),\, z\in E.
\]
Outside $B_2(z)$, we have $\Phi_z=A_z$. Since $\Phi_z:\mathbb R^D\to \mathbb R^D$ is one to one and onto, it follows that
$\Phi_z: B_2(z)\to A_z(B_2(z))$ is one to one and onto. Consequently, since $\Phi=\Phi_z$ on $B_2(z)$, we have that:
\begin{itemize}
\item $\Phi: B_2(z)\to A_z (B_2(z))$ is one to one and onto for each $z\in E$.
\end{itemize}
Next, recall that $\Phi_z^*=A_z$ on $B_2(z)$, in particular
\[
\Phi_z^*: B_2(z)\to A_z( B_2(z))
\]
is one to one and onto. Also, $\Phi_z^*:\mathbb R^D\to \mathbb R^D$ is one to one and onto and $\Phi_z^*=A_z^*$ outside $B_4(z)$ so it follows that
\[
\Phi_z^*:B_4(z)\to A_z^*(B_4(z))
\]
is one to one and onto. Consequently
\[
\Phi_z^*:B_4(z)\setminus B_2(z)\to A_z^{*}(B_4(z))\setminus A_z (B_2(z))
\]
is one to one and onto. Since $\Phi=\Phi^{*}$ on $B_4(z)\setminus B_2(z)$, we conclude that
\begin{itemize}
\item $\Phi: B_4(z)\setminus B_2(z)\to A_z^* (B_4(z))\setminus A_z (B_2(z))$ is one to one and onto for $z\in E$.
\end{itemize}

Next, recall that $\Psi:\mathbb R^D\to\mathbb R^D$ is one to one and onto and that $\Psi=A_z^*$ on $B_4(z)$ for each
$z\in E$. Hence,
\[
\Psi:\mathbb R^D\setminus \cup_{z\in E}B_4(z)\to \mathbb R^D\setminus \cup_{z\in E}A_z^* (B_4(z))
\]
is one to one and onto. Since $\Phi=\Psi$ on $\mathbb R^D\setminus \cup_{z\in E}B_4(z)$, we conclude that
\begin{itemize}
\item
\[
\Phi:\mathbb R^D\setminus \cup_{z\in E}B_4(z) \to \mathbb R^D\setminus \cup_{z\in E}A_z^*(B_4(z))
\]
is one to one and onto.
\end{itemize}
Recall that $B_2(z)\subset B_4(z)$ for each $z\in E$ and that the balls $B_4(z), z\in E$ are pairwise disjoint. Hence the following sets constitute a partition of $\mathbb R^D$:

\begin{itemize}
\item $B_2(z)$ (all $z\in E$); $B_4(z)\setminus B_2(z)$ (all $z\in E$); $\mathbb R^D\setminus \cup_{z\in E}B_4(z)$.
\end{itemize}

Moreover, we recall that $A_z,A_z^*$ are Euclidean motions, $B_2(z), B_4(z)$ are balls centered at $z$ with radii $\exp(-3/\varepsilon)\tau$ and $\exp(-1/\varepsilon)\tau$ respectively and
\[
|A_z(z)-A_z^*(z)|\leq C\varepsilon\exp(-4/\varepsilon)\tau.
\]
It follows that $A_z (B_2(z))\subset A_z^* (B_4(z))$ for $z\in E$. Moreover, $A_z^*=\phi(z)$ for $z\in E$. For $z,z'\in E$ distinct, we have
\[
|\phi(z)-\phi(z')|\geq 1/2|z-z'|\geq 1/2\tau.
\]
Since, $A_z^* (B_4(z))$ is a ball of radius $\exp(-1/\varepsilon)\tau$ centered at $\phi(z)$ for each $z\in E$, it follows that the balls $A_z^* (B_4(z))\, (z\in E)$ are pairwise disjoint. Therefore the 
following sets constitute a partition of $\mathbb R^D$:
\begin{itemize}
\item $A_z (B_2(z))(z\in E)$,\, $A_z* (B_4(z))\setminus A_z (B_2(z))(z\in E)$,\, $\mathbb R^D\setminus \cup_{z\in E}A_z^*( B_4(z))$.
\end{itemize}

In view of the 5 bullet points regarding the partitions of $\mathbb R^D$ and the bijective character of $\Phi$ restricted appropriately, we conclude that $\Phi:\mathbb R^D\to\mathbb R^D$ is one to one and onto. 
The proof of the Gluing Theorem is complete. $\Box$.

\section{Hierarchical clusterings of finite subsets of $\mathbb R^D$.}
\setcounter{equation}{0}
We are almost ready for the proofs of Theorem~\ref{t:Theorem2a} and Theorem~\ref{t:Theorem2b}. We need one more piece of machinery
namely a technique developed by Fefferman in \cite{F1,F2,F3,F4,F5} which involves the combinatorics of hierarchical clusterings of finite subsets of $\mathbb R^D$. This is also of independent interest.
\subsection{Hierarchical clusterings of finite subsets of $\mathbb R^D$.}
Our main result in this section is the following hierarchical clustering result. 

\begin{lem}
Let $S\subset \mathbb R^D$ with $2\leq {\rm card}(S)\leq K$. Here $K$ is given as in (~\ref{e:specialk}). Let $\varepsilon>0$. Then there exists $\tau$ satisfying
\[\exp(-C_K/\varepsilon){\rm diam}(S)\leq \tau\leq \exp(-1/\varepsilon){\rm diam}(S)\]
and a partition of $S$ into subsets $S_{\nu}(\nu=1,...,\nu_{(max)})$ with the following properties:
\begin{itemize}
\item ${\rm card}(S_{\nu})\leq K-1,\, \forall \nu.$
\item ${\rm diam}(S_{\nu})\leq \exp(-5/\varepsilon)\tau,\, \forall \nu.$
\item ${\rm dist}(S_{\nu}, S_{\nu'})\geq \tau,\, \forall \nu.$
\label{l:clustering}
\end{itemize}
\end{lem}

{\bf Proof:}\, We define and equivalence relation on $S$ as follows. Define a relation $\sim$ on $S$ by saying that $x\sim x'$, for $x,x'\in S$ if and only if 
$|x-x'|\leq \exp(-5/\varepsilon)\tau$ for a $\tau>0$ to be chosen in a moment. By the pigeonhole principle, we may choose and fix a $\tau$ satisfying 
\[\exp(-C_K/\varepsilon){\rm diam}(S)\leq \tau\leq \exp(-1/\varepsilon){\rm diam}(S)\]
so that $\sim$ is an equivalence relation for such fixed $\tau$. Then the equivalence classes of $\sim$ partitions $S$
into the sets with the properties as required. $\Box$.

The same technique as above yields,  see (\cite{DF4}):

\begin{lem}
Let $k\geq 2$ be a positive integer and let $0<\eta\leq 1/10$. Let $S\subset \mathbb R^n$ be a set consisting of $k$
distinct points. Then, we can partition $S$ into sets $E_1,E_2,...E_{{\mu}_{\rm max}}$ and we can find a positive integer $l$
$(10\leq l\leq 100+\binom{k}{2})$ such that the following hold:
\beq
{\rm diam}(S_{\mu})\leq \eta^{l}{\rm diam}(S)
\eeq
for each $\mu$ and
\beq
{\rm dist}(S_{\mu},S_{\mu'})\geq \eta^{l-1}{\rm diam}(S),\, {\rm for}\, \mu\neq \mu'.
\eeq
\label{l:lemmathree}
\end{lem}

\section{Proof of Theorem~\ref{t:Theorem2a} and Theorem~\ref{t:Theorem2b}.}
\setcounter{equation}{0}
We begin with the
\medskip

{\bf Proof of Theorem~\ref{t:Theorem2a}}:\, We use induction on $K$. If $K=1$, the theorem holds trivially. For the induction step
we will fix $K\geq 2$ and assume that our result holds for $K-1$. We now establish the theorem for the given
$K$. Thus, we are making the following inductive assumptions. For suitable constants $c_{\rm old}$, $C_{\rm old}', C_{\rm old}''$ depending only on
$D,K$ the following holds: Inductive hypothesis: Suppose that $0<\varepsilon<c_{\rm old}$, define
$\eta_{\rm old}=\exp(-C_{\rm old}'/\varepsilon)$ and $\delta_{\rm old}=\exp(-C_{\rm old}''/\varepsilon)$. Let $\phi^*:S^*\to \mathbb R^D$ with $S^*\subset \mathbb R^D$ and ${\rm card}(S^*)\leq K-1$. Suppose \[
(1+\delta_{\rm old})^{-1}|x-y|\leq |\phi^*(x)-\phi^*(y)|\leq (1+\delta_{\rm old})|x-y|,\, x,y\in S. \]
Then the following holds: If $\phi^*$ has no negative $\eta_{\rm old}$ block, then there exists a proper $\varepsilon$-distorted diffeomorphism $\Phi^*:\mathbb R^D\to \mathbb R^D$ such that $\phi^*=\Phi^*$ on $S$ and $\Phi^*$ agrees with a proper Euclidean motion on
\[
\left\{x\in \mathbb R^D:{\rm dist}(x,S^*)\geq 10^4{\rm diam}(S^*)\right\}.
\]
Now let $L$, $L'$, $L''$ be positive constants to be fixed later. (Eventually we will let them depend on $D$ and $K$ but not yet). Now suppose that
\begin{itemize}
\item[(1)] $0<\varepsilon<L$.
\item[(2)] We set $\eta=\exp(-L'/\varepsilon)$
\item[(3)] and we set $\delta=\exp(-L''/\varepsilon)$.
\item[(4)] Let $\phi:S\to \mathbb R^D$
\item[(5)] where $S\subset \mathbb R^D$
\item[(6)] ${\rm card}(S)=K$ and
\item[(7)] $(1+\delta)^{-1}|x-y|\leq |\phi(x)-\phi(y)|\leq (1+\delta)|x-y|,\, x,y\in S.$
\item[(8)] Suppose that $\phi$ has a negative $\eta$ block.
\end{itemize}
We will construct a proper $C\varepsilon$-distorted diffeomorphism
$\Phi$ that agrees with $\phi$ on $S$ and with a proper Euclidean motion away from $S$. To do, so we first apply the clustering lemma,
Lemma~\ref{l:clustering}. Recall that ${\rm card}(S)=K\geq 2$ so the clustering lemma applies. Let $\tau$ and $S_{\nu}(\nu=1,...,\nu_{\rm max})$ be as in the clustering lemma. Thus,
\begin{itemize}
\item[(9)] $S$ is the disjoint union of $S_{\nu}(\nu=1,...,\nu_{\rm max})$.
\item[(10)] ${\rm card}(S_{\nu})\leq K-1$ for each $\nu(\nu=1,...,\nu_{\rm max})$.
\item[(11)] ${\rm diam}S_{\nu}\leq \exp(-5/\varepsilon)\tau$ for each $\nu(\nu=1,...,\nu_{\rm max})$.
\item[(12)] ${\rm dist}(S_{\nu}, S_{\nu'})\geq \tau{\rm diam}(S)$, for $\nu\neq \nu'$,
for each $\nu,\nu'(\nu,\nu'=1,...,\nu_{\rm max})$.
\item[(13)] $\exp(-C_K/\varepsilon){\rm diam}(S)\leq \tau\leq \exp(-1/\varepsilon){\rm diam}(S).$
\item[(14)] Assuming that $L'>C'_{\rm old}$ and $L''>C''_{\rm old}$, we see that $\eta<\eta_{\rm old}$ and $\delta<\delta_{\rm old}$. Hence by (7) and (8) we have:
\item[(15)] $\phi|S_\nu$ does not have an $\eta_{\rm old}$ block and
\[
(1+\delta_{\rm old})^{-1}|x-y|\leq |\phi(x)-\phi(y)|\leq (1+\delta_{\rm old})|x-y|,\, x,y\in S_{\nu} \]
Consequently (10) and the induction hypothesis
\item[(16)] produce a proper $\varepsilon$-distorted diffeomorphism $\Phi_{\nu}:\mathbb R^D\to \mathbb R^D$ such that
\item[(17)] $\Phi_{\nu}=\phi$ on $S_{\nu}$ and
\item[(18)] $\Phi_{\nu}=A_{\nu}$ on
$\left\{x\in \mathbb R^D: {\rm dist}(x, S_{\nu})\geq 10^4{\rm diam}(S_{\nu})\right\}$
where $A_{\nu}$ is a proper Euclidean motion.
\end{itemize}
Next, for each $\nu$ $(1\leq \nu\leq \nu_{{\rm max}})$, we pick a representative $y_{\nu}\in S_{\nu}$. Define
\begin{itemize}
\item[(19)] $E=\left\{y_{\nu}:\, 1\leq \nu\leq \nu_{{\rm max}}\right\}$.
\item[(20)] Thus $E\subset \mathbb R^D$, \, $2\leq {\rm card}(E)\leq K$,
\item[(21)] $\frac{1}{2}{\rm diam}(S)\leq {\rm diam}(E)\leq {\rm diam}(S)$ and by (12) and (13),
\item[(22)]
\[
|x-y|\geq \tau \geq \exp(-C_K/\varepsilon){\rm diam}(S)
\]
for $x,y\in E$ distinct.
\end{itemize}

We prepare to apply a rescalled version of Theorem~\ref{t:lemmareflectionex4}. For easier reading, let us note the assumptions and conclusions with the same notation there as we will need to verify and use them here.

\begin{itemize}
\item {\bf Assumptions on E}.
\item[(23)] ${\rm card}(E)\leq K$
\item[(24)] $|x-y|\geq \tau$ for $x,y\in E$ distinct.
\item {\bf Assumptions on $\phi$.}
\item[(25)] \[
(1+\delta)^{-1}|x-y|\leq |\phi(x)-\phi(y)|\leq (1+\delta)|x-y|,\ , x,y\in E.
\]
\item[(26)] $\phi$ has no negative $\eta$ blocks.
\item {\bf Assumptions on the parameters.}
\item[(27)] $0<\eta<c\varepsilon\tau/{\rm diam}(E)$ for small enough $c$
\item[(28)] $C_K\delta^{1/\rho_K}\tau^{-1}{\rm diam}(E)\leq {\rm min}(\varepsilon,\eta^D)$ for large enough $C_K, \rho_K$ depending only on $K$ and $D$.
\item {\bf Conclusion.}
\item[(28a)] There exists a proper $C\varepsilon$-distorted diffeomorphism $\Psi:\mathbb R^D\to \mathbb R^D$ with the following properties:
\item[(29)] $\Psi=\phi$ on $E$.
\item[(30)] $\Psi$ agrees with a proper Euclidean motion on
\[
\left\{x\in \mathbb R^D:\, {\rm dist}(x,E)\geq 1000{\rm diam}(E)\right\}.
\]
\item[(31)] For each $z\in E$, $\Phi$ agrees with a proper Euclidean motion on $B(z,\tau/1000)$.
\end{itemize}

Let us check that our present $\phi:E\to \mathbb R^D$, $\delta, \varepsilon, \eta, \tau$ satisfy the hypotheses of
Theorem~\ref{t:lemmareflectionex4}. In fact: Hypothesis (23) is (20). Hypothesis (24) is (22),
Hypothesis (25) is immediate from (7).
Hypothesis (26) is immediate from (8).

Let us check hypotheses (27) and (28). From
(13) and (21) we have
\begin{itemize}
\item[(32)]
\[
\exp(-C_K/\varepsilon)\leq \tau/{\rm diam}(E). \]
Hence (27) and (28) will follow if we can show that the following two things:
\item[(33)] $0<\eta<c\exp(-C_K/\varepsilon)$ for small enough $c$.
\item[(34)] $C_K\delta^{1/\rho_k}\exp(C_K/\varepsilon)\leq {\rm min}(\varepsilon,\eta^D)$.
However, we now recall that $\delta$ and $\eta$ are defined by (2) and (3). Thus (33) holds provided
\item[(35)] $L<c_K$ for small enough $c_K$ and $L'>C_K$ for large enough $C_K$.
\item[(36)] Similarly (34) holds provided $L<c_K$ for small enough $c_K$ and $1/\rho_K L''-C_K\geq {\rm max}(1, DL')$.
Assuming we can choose $L,L',L''$ as we wish, we have (33) and (34) hence also (27) and (28). This completes our verification of the hypothesis of Theorem~\ref{t:lemmareflectionex4} for our present $\Phi$ and $E$. Applying Theorem~\ref{t:lemmareflectionex4}, we now obtain a proper $C\varepsilon$-distorted diffeomorphism $\Phi:\mathbb R^D\to \mathbb R^D$ satisfying (28a-31). For each $z\in E$, we now define a proper $\varepsilon$-distorted diffeomorphism $\Phi_z$ by setting:
\item[(37)] $\Phi_z=\Phi_{\nu}$ if $z=y_{\nu}$. (Recall (16), (19) and note that the $y_{\nu}$, $1\leq \nu\leq \nu_{{\rm max}}$ are distinct). From (17), (18), (37) we have the following:
\item[(38)] $\Phi_z=\phi$ on $S_{\nu}$ if $z=y_{\nu}$. In particular,
\item[(39)] $\Phi_z(z)=\phi(z)$ for each $z\in E$. Also
\item[(40)] $\Phi_z=A_z$ (a proper Euclidean motion) outside $B(z, 10^5{\rm diam}(S_{\nu}))$ if $z=y_{\nu}$. Recalling
(11), we see that
\item[(41)] $\Phi_z=A_z$ (a proper Euclidean motion) outside $B(z, 10^5\exp(-5/\varepsilon)\tau)$. We prepare to apply Theorem~\ref{t:lemmagl}, the Gluing Lemma to the present $\phi$, $E$, $\Phi_z(z\in E)$, $\Psi$, $\varepsilon$ and $\tau$. Let us check the hypotheses of the Gluing Lemma. We have $\phi:E\to \mathbb R^D$ and
$1/2|x-y|\leq |\phi(x)-\phi(y)|\leq 2|x-y|$ for $x,y\in E$ thanks to (7) provided
\item[(41a)] $L\leq 1$ and $L''\geq 10$. See also (3). Also $|x-y|\geq \tau$ for $x,y\in E$ distinct, see (22). Moreover, for each $z\in E$, $\Phi_z$ is a proper $\varepsilon$ distorted diffeomophism (see (16) and (37)). For each $z\in E$, we have $\Phi_z(z)=\phi(z)$ by (39) and
$\Phi_z=A_z$ (a proper Euclidean motion) outside $B_1(z)=B(z,\exp(-4/\varepsilon)\tau)$, see (41). Here, we assume that,
\item[(42)] $L\leq c_k$ for a small enough $c_k$.
Next, recall that $\Psi$ satisfies (28a-31). Then $\Psi$ is a $C\varepsilon$-distorted diffeomorphism, $\Psi=\Phi$ on $E$ and for $z\in E$, $\Psi$ agrees with a proper Euclidean motion $A_z^*$ on $B(z,\frac{\tau}{1000})$, hence on
$B_4(z)=B(z, \exp(-1/\varepsilon)\tau)$. Here, again, we assume that $L$ satisfies (42). This completes the verification of the hypotheses of the Gluing Lemma. Applying that lemma, we obtain:
\item[(43)] a $C' \varepsilon$-distorted diffeomorphism $\Phi:\mathbb R^D\to \mathbb R^D$ such that
\item[(44)] $\Phi=\Phi_z$ on $B_2(z)=B(z,\exp(-3/\varepsilon)\tau)$, for each $z\in E$ and
\item[(45)] $\Phi=\Psi$ outside $\cup_{z\in E}B_3(z)=\cup_{z\in E}B(z,\, \exp(-2/\varepsilon)\tau).$
Since $\Psi$ is proper, we know that
\item[(46)] $\Phi$ is proper.
Let $z=y_{\mu}\in E$. Then (11) shows that $S_{\mu}\subset B(z,\, \exp(-5/\varepsilon)\tau)$ and therefore (44) yields
$\Phi=\Phi_z$ on $S_{\mu}$ for $z=y_{\mu}$. Together, with (38), this yields $\Psi=\phi$ on $S_{\mu}$ for each $\mu (1\leq \mu\leq \mu_{\rm max})$. Since the $S_{\mu}(1\leq \mu\leq \mu_{\rm max})$ form a partition of $S$, we conclude that
\item[(47)] $\Phi=\phi$ on $S$. Moreover, suppose that
\[
{\rm dist}(x, S)\geq 10^4{\rm diam}(S).
\]
Then $x$ does not belong to $B(z, \, \exp(-2/\varepsilon)\tau)$ for any $z\in E$ as we see from (13). (Recall that $E\subset S$).
Consequently (45) yields $\Phi(x)=\Psi(x)$ and therefore (30) tells us that $\Psi(x)=A_{\infty}(x)$.  Since
${\rm dist}(x, S)\geq 10^4{\rm diam}(S)$, we have
\[
{\rm dist}(x, E)\geq {\rm dist}(x, S)\geq 10^4{\rm diam}(S)\geq 10^3{\rm diam E}.
\]
Hence (30) applies. Thus,
\item[(48)] $\Phi$ agrees with a proper Euclidean motion $A_{\infty}^*$ on
\[
\left\{x\in \mathbb R^D:\, {\rm dist}(x, S)\geq 10^4{\rm diam}(S)\right\}.
\]
Collecting our results (43), (46), (47), (48), we have the following:
\item[(49)] There exists a proper $C’\varepsilon$-distorted diffeomorphism $\Phi$ such that $\Phi=\phi$ on $S$ and $\Phi$ agrees with a proper Euclidean motion on
\[
\left\{x\in \mathbb R^D:\, {\rm dist}(x, S)\geq 10^4{\rm diam}(S)\right\}.
\]
We have established (49) assuming that the small constant $L$ and the large constants $L'$, $L''$ satisfy the conditions (14), (35), (36), (41a), (42). By picking $L$ first, $L'$ second and $L''$ third, we can satisfy all those conditions with $L=c_K$, $L'=C'_K$,
$L''=C''_K$. With these $L, L', L''$ we have shown that (1)-(8) together imply that (49) holds. Thus, we have proven the following:
\item[(50)] For suitable constants $C$, $C_K$, $C_K'$, $C_K''$ depending only on $D$ and $K$ the following holds: Suppose that
$0<\varepsilon<c_k$. Set $\eta=\exp(-C_K'/\varepsilon)$ and $\delta=\exp(-C_K''/\varepsilon)$. Let $\phi:S\to \mathbb R^D$
with ${\rm card}(S)=K$, $S\subset \mathbb R^D$. Assume that
\[
(1+\delta)^{-1}|x-y|\leq |\phi(x)-\phi(y)|\leq (1+\delta)|x-y|, \, x, y\in S.
\]
Then if $\phi$ has no negative $\eta$ block, then there exists a proper $C\varepsilon$-distorted diffeomorphism
$\Phi:\mathbb R^D\to \mathbb R^D$ such that $\phi=\Phi$ on $S$ and $\Phi$ agrees with a proper Euclidean motion on
\[
\left\{x\in \mathbb R^D:\, {\rm dist}(x, S)\geq 10^4{\rm diam}(S)\right\}.
\]
Taking $\varepsilon$ to be $\varepsilon/C$, we thus deduce:
\item[(51)] For suitable constants $C_{{\rm new}}$, $C_{{\rm new}}'$, $C_{{\rm new}}''$ depending only on $D$ and $K$ the following holds: Suppose that
$0<\varepsilon<c_{\rm new}$. Set $\eta=\exp(-C_{{\rm new}}'/\varepsilon)$ and $\delta=\exp(-C_{{\rm new}}''/\varepsilon)$. Let $\phi:S\to \mathbb R^D$
with ${\rm card}(S)=K$, $S\subset \mathbb R^D$. Assume that
\[
(1+\delta)^{-1}|x-y|\leq |\phi(x)-\phi(y)|\leq (1+\delta)|x-y|, \, x, y\in S.
\]
Then if $\phi$ has no negative $\eta$ block, then there exists a proper $\varepsilon$-distorted diffeomorphism
$\Phi:\mathbb R^D\to \mathbb R^D$ such that $\phi=\Phi$ on $S$ and $\Phi$ agrees with a proper Euclidean motion on
\[
\left\{x\in \mathbb R^D:\, {\rm dist}(x, S)\geq 10^4{\rm diam}(S)\right\}.
\]
\end{itemize}
Thats almost Theorem~\ref{t:Theorem2a} except we are assuming ${\rm card}(S)=K$ rather than ${\rm card}(S)\leq K$. Therefore we proceed as follows: We have our result (50) and we have an inductive hypothesis. We now take
$C'={\rm max}(C_{{\rm old}}', C_{{\rm new}}')$, $C''={\rm old}(C_{{\rm max}}'', C_{{\rm new}}'')$ and
$c'={\rm min}(c_{{\rm old}}', c_{{\rm new}}')$. These constants are determined by $D$ and $K$. We now refer to
$\eta_{{\rm old}}$, $\eta_{{\rm new}}$, $\eta$, $\delta_{{\rm old}}$, $\delta_{{\rm new}}$, $\delta$ to denote
$\exp(-C_{{\rm old}}'/\varepsilon)$, $\exp(-C_{{\rm new}}'/\varepsilon)$, $\exp(-C'/\varepsilon)$, $\exp(-C_{{\rm old}}''/\varepsilon)$, $\exp(-C_{{\rm new}}''/\varepsilon)$, $\exp(-C''/\varepsilon)$ respectively.

Note that $\delta\leq \delta_{{\rm old}}$, $\delta\leq \delta_{{\rm new}}$, $\eta\leq \eta_{{\rm old}}$ and $\eta\leq \eta_{{\rm new}}$. Also
if $0<\varepsilon<c$, then $0<\varepsilon<c_{{\rm old}}$ and $0<\varepsilon<c_{{\rm new}}$. If
\[
(1+\delta)^{-1}|x-y|\leq |\phi(x)-\phi(y)|\leq (1+\delta)|x-y|, x,y\in S
\]
then the same holds for $\delta_{{\rm old}}$ and $\delta_{{\rm new}}$. Also if $\phi$ has no negative $\eta$-block, then it has no negative
$\eta_{{\rm old}}$-block and it has no negative $\eta_{{\rm new}}$-block. Consequently by using (51) and the induction hypothesis, we have proved Theorem~\ref{t:Theorem2a}. $\Box$
\medskip

We now give
\medskip

{\bf The Proof of Theorem~\ref{t:Theorem2b}}:\, To see this, we simply observe that increasing $C_K''$ in the Theorem~\ref{t:Theorem2a} above merely weakens the result so we may increase $C_K''$ and achieve that $0<\delta<c\eta^D$ for small enough $c$ and $0<\delta<\varepsilon$. The desired result then follows by using Corollary~\ref{c:cextensionblock1}. $\Box$

\section{ Proofs of Theorem~\ref{t:Theorem3} and Theorem~\ref{t:Theorem4}.}
\setcounter{equation}{0}

It remains to give the proofs of Theorem~\ref{t:Theorem3} and Theorem~\ref{t:Theorem4}.
\medskip

{\bf The Proof of Theorem~\ref{t:Theorem3}}:\, Pick $c_K$, $C_K$, $C_K''$ as in Theorem~\ref{t:Theorem2a} and Theorem~\ref{t:Theorem2b}. Let $\delta$ and $\eta$ be as in Theorem~\ref{t:Theorem2a} and let us take $S_0=\left\{x,y\right\}$. Then we see that
\[
(1+\delta)|x-y|\leq |\phi(x)-\phi(y)|\leq (1+\delta)|x-y|, x, y\in S.
\]
Now, if $\phi$ has no negative $\eta$-block, then by Theorem~\ref{t:Theorem2a}, $\Phi$ exists with the properties claimed. Similarly,
if $\phi$ has no positive $\eta$-block, then applying Theorem~\ref{t:Theorem2a} to the map $\phi o({\rm reflection})$ we obtain the $\Phi$ we need with the properties claimed. Suppose that $\phi$ has a positive $\eta$-block $(x_0,...x_D)$ and a negative $\eta$-block $(y_0,...y_D)$.
Then by Theorem~\ref{t:Theorem2b}, $\phi|_{\left\{x_0,...,x_D, y_0,...,y_D\right\}}$ cannot be extended to a $\delta$ distorted diffeomorphism
$\Phi:\mathbb R^D\to \mathbb R^D$. Indeed, the $\eta$-block $(x_0,...,x_D)$ forces any such $\Phi$ to be proper while the $\eta$-block $(y_0,...y_D)$ forces $\Phi$ to be improper. Since ${\rm card}\left\{x_0,...,x_D,y_0,...y_D\right\}\leq 2(D+1)$, the proof of Theorem~\ref{t:Theorem3} is complete. $\Box$
\medskip

{\bf The Proof of Theorem~\ref{t:Theorem4}}:\, Take $k\leq D+1$ and apply Theorem~\ref{t:Theorem2a}. Let $\eta$ and $\delta$ be determined by
$\varepsilon$ as in Theorem~\ref{t:Theorem2a}. If $\phi$ has no negative $\eta$-block then applying Theorem~\ref{t:Theorem2a} to $\phi$ or
$\phi o({\rm reflection})$, we see that $\phi$ extends to a $\varepsilon$-distorted diffeomorphism $\Phi:\mathbb R^D\to \mathbb R^D$.
However since $({\rm card})(S)\leq D+1$, the only possible (negative or positive) $\eta$-block for $\phi$ is all of $S$. Thus either $\phi$ has a negative $\eta$-block or it has no positive $\eta$-block. $\Box$

\noindent
{\bf Acknowledgment:}\, Support from the Department of Mathematics at Princeton, the National Science Foundation, The American Mathematical Society, The Department of Defense and the University of the Witwatersrand are gratefully acknowledged.

\end{document}